\documentclass[reqno]{amsart}

\usepackage{lmodern}
\usepackage{libertine}
\usepackage{libertinust1math}
\usepackage[utf8]{inputenc}
\usepackage{geometry}[margin=1in]
\usepackage{amsmath,amsthm,amssymb,amsfonts, color,dsfont,enumitem}
\usepackage[hidelinks]{hyperref}
\usepackage[normalem]{ulem}
\usepackage{tikz}
\usepackage{multicol}
\usepackage{float}
\usepackage{siunitx,booktabs,array}

\usetikzlibrary{braids}

\newcommand{\CC}{\mathcal{C}}
\newcommand{\DD}{\mathcal{D}}

\newcommand{\ZZ}{\mathcal{Z}}

\title{Zesting produces modular isotopes and explains their topological invariants}
\author{Colleen Delaney, Sung Kim, and Julia Plavnik}

\address{Department of Mathematics, 
Department of Physics and Astronomy, Purdue University, West Lafayette, IN 47907 USA}
\email{colleend@purdue.edu}

\address{Department of Mathematics, University of Southern California, Los Angeles, CA 90089 USA}
\email{skim2261@usc.edu}

\address{Department of Mathematics, Indiana University, Bloomington, IN 47405 USA}
\email{jplavnik@iu.edu}

\makeatletter
\@namedef{subjclassname@2020}{
  \textup{2020} Mathematics Subject Classification}
\makeatother

\newtheorem{defn}[]{Definition}
\newtheorem{lem}[]{Lemma}
\newtheorem{prop}[]{Proposition}
\newtheorem{quest}[]{Question}
\newtheorem{rmk}[]{Remark}
\newtheorem{thm}[]{Theorem}

\theoremstyle{definition}
\newtheorem{examp}[]{Example}

\DeclareMathOperator{\Aut}{Aut}

\DeclareMathOperator{\FPdim}{FPdim}
\DeclareMathOperator{\id}{id}
\DeclareMathOperator{\Id}{Id}
\DeclareMathOperator{\Hom}{Hom}
\DeclareMathOperator{\Inv}{Inv}
\DeclareMathOperator{\Irr}{Irr}
\DeclareMathOperator{\rank}{rank}
\DeclareMathOperator{\Rep}{Rep}
\DeclareMathOperator{\Tr}{Tr}
\DeclareMathOperator{\Vect}{Vec}

\newcommand{\parcen}[1]{\phantom{(}#1\phantom{)}}
\newcommand{\br}{to [out=90,in=-90]}

\newcommand{\twobox}[5]{
\begin{scope}[line width=1]
\draw (0,0) node[below] {$#1$} --(0,2) node[above] {$#2$};
\draw (1,0) node[below] {$#3$} --(1,2) node[above] {$#4$};
\draw[fill=white] (-.35,.5) rectangle node {$#5$} (1.35,1.5);
\end{scope}
}

\newcommand{\cocycm}[5]{
  \foreach \x in {0,1.5}{
  \draw(\x,0)--(\x,2);
  }
  \draw[fill=white] (-.25,.75) rectangle (1.75,1.25); 
  \draw (0,0) node[below] {$#4(#1,#2#3)$};
  \draw (1.5,0) node[below] {$#4(#2,#3)$};
  \draw (0,2) node[above] {$#4(#1,#2)$};
  \draw (1.5,2) node[above] {$#4(#1#2,#3)$};
  \draw (.75,1) node {$#5_{#1,#2,#3}$};
  }

  \newcommand{\kink}[2]{
\begin{scope}[line width=1]
\draw (0,0) node[below] {$#1$}-- (0,1) to [out=90,in=180] (.25,1.25) to [out=0,in=90] (.5,1) to [out=-90,in=0] (.25,.75);
\draw[white,line width=8] (.25,.75) to [out=180,in=-90] (0,1)--(0,2);
\draw (.25,.75) to [out=180,in=-90] (0,1)--(0,2) node[above] {$#2$};
\end{scope}
}

  \newcommand{\onebox}[3]{
\begin{scope}[line width=1]
\draw (0,0) node[below] {$#1$} --(0,2) node[above] {$#2$};
\draw[fill=white] (-.33,.66) rectangle node {#3} (.33,1.33);
\end{scope}
}

  \newcommand{\onecirc}[3]{
\begin{scope}[line width=1]
\draw (0,0) node[below] {$#1$} --(0,2) node[above] {$#2$};
\draw[fill=white] (0,1) circle (.45) node {$#3$};
\end{scope}
}

\newcommand{\seqassoc}[3]{\draw (1,.5) \br (2,2);
\draw[white, line width=10] (2,.5) \br (1,2);
\draw (2,.5) \br (1,2);
\begin{scope}[yshift=3.5cm]
\draw (2,.5) \br (1,2);
\draw[white, line width=10] (1,.5) \br (2,2);
\draw (1,.5) \br (2,2);
\end{scope}
\begin{scope}[yshift=2cm]
\twobox{}{}{}{}{\nu_{#1,#2,#3}}
\end{scope}
\begin{scope}[xshift=2cm, yshift=2cm]
\twobox{}{}{}{}{\nu'_{#1,#2,#3}}
\end{scope}
\foreach \x in {0,3}
{\draw (\x, .5)--(\x,2);
\draw (\x,4) -- (\x,5.5);}}

\newcommand{\seqassocmodbottom}[3]{
\draw[white, line width=10] (2,.5) \br (1,2);
\draw (2,.5) \br (1,2);
\begin{scope}[yshift=3.5cm]
\draw (2,.5) \br (1,2);
\draw[white, line width=10] (1,.5) \br (2,2);
\draw (1,.5) \br (2,2);
\end{scope}
\begin{scope}[yshift=2cm]
\twobox{}{}{}{}{\nu}
\end{scope}
\begin{scope}[xshift=2cm, yshift=2cm]
\twobox{}{}{}{}{\nu'}
\end{scope}
\foreach \x in {0,3}
{\draw (\x, .5)--(\x,2);
\draw (\x,4) -- (\x,5.5);}}

\newcommand{\seqassocmodtop}[3]{\draw (1,.5) \br (2,2);
\draw[white, line width=10] (2,.5) \br (1,2);
\draw (2,.5) \br (1,2);
\begin{scope}[yshift=3.5cm]
\draw[white, line width=10] (1,.5) \br (2,2);
\draw (1,.5) \br (2,2);
\end{scope}
\begin{scope}[yshift=2cm]
\twobox{}{}{}{}{\nu}
\end{scope}
\begin{scope}[xshift=2cm, yshift=2cm]
\twobox{}{}{}{}{\nu'}
\end{scope}
\foreach \x in {0,3}
{\draw (\x, .5)--(\x,2);
\draw (\x,4) -- (\x,5.5);}}

\newcommand{\zestedbraid}[2]{
\begin{scope}[line width=1]
    \foreach \x in {0,3,4}{
        \draw(\x,#1)--(\x,#1 + 1);
    }

    \draw[fill=white] (3,#1 + .5) circle (.33) node {\small $#2$};
      \draw (1,#1) to [out=90, in=-90] (2,#1 + 1);
      \draw[white, line width=10] (2,#1)  to [out=90, in=-90] (1,#1 + 1);
    \draw (2,#1)  to [out=90, in=-90] (1,#1 + 1);
\end{scope}
}

\newcommand{\zestedbraidminv}[2]{
\begin{scope}[line width=1]
    \foreach \x in {2,3,4}{
        \draw(\x,#1)--(\x,#1 + 1);
    }
    \draw[fill=white] (2,#1 + .5) circle (.33) node {\small $#2$};
    \draw (2.2, #1 + .85) node[right] {$\scriptstyle{-1}$};
    \draw (1,#1)  to [out=90, in=-90] (0,#1 + 1);
      \draw[white, line width=10] (0,#1)  to [out=90, in=-90] (1,#1 + 1);
         \draw (0,#1) to [out=90, in=-90] (1,#1 + 1);
\end{scope}
}

\newcommand{\zestedassoc}[4]{
\begin{scope}[line width=1]
    \foreach \x in {0,1,4}{
        \draw(\x,#4)--(\x,#4 + 2);
    }
    \foreach \x in {2,3}{
        \draw(\x,#4)--(\x,#4 + 1);
    }
    \draw(3,#4+1)  to [out=90, in=-90] (2,#4 + 2);
    \draw[white, line width=10] (2,#4+1)  to [out=90, in=-90] (3,#4 + 2);
    \draw (2,#4+1) to [out=90, in=-90] (3,#4 + 2);
    \draw[fill=white] (2.75,#4 + .2) rectangle node {$\nu_{#1,#2,#3}$} (4.25,#4 + .8);
    \end{scope}
}

\newcommand{\zestedassocinv}[4]{
\begin{scope}[line width=1]
    \foreach \x in {0,1,4}{
        \draw(\x,#4)--(\x,#4 + 2);
    }
    \foreach \x in {2,3}{
        \draw(\x,#4+1)--(\x,#4 + 2);
    }
     \draw (2,#4) to [out=90, in=-90] (3,#4 + 1);
      \draw[white, line width=10] (3,#4)  to [out=90, in=-90] (2,#4 + 1);
    \draw(3,#4)  to [out=90, in=-90] (2,#4 + 1);
    \draw[fill=white] (2.75,#4 + 1.2) rectangle node {$\nu_{#1,#2,#3}$} (4.25,#4 + 1.8);
 \draw (4.25, #4 + 1.9) node[right] {$\scriptstyle{-1}$};
 \end{scope}
}

\newcommand{\zestedseqcochain}[3]{
\begin{scope}[line width=1, xscale=1, yscale=.5]
\draw (1,.5) \br (2,2);
\draw[white, line width=10] (2,.5) \br (1,2);
\draw (2,.5) \br (1,2);
\begin{scope}[yshift=3.5cm]
\draw (2,.5) \br (1,2);
\draw[white, line width=10] (1,.5) \br (2,2);
\draw (1,.5) \br (2,2);
\end{scope}
\begin{scope}[yshift=2cm]
\twobox{}{}{}{}{\nu_{#1,#2,#3}}
\end{scope}
\begin{scope}[xshift=2cm, yshift=2cm]
\twobox{}{}{}{}{\nu'_{#1,#2,#3}}
\end{scope}
\foreach \x in {0,3}
{\draw (\x, .5)--(\x,2);
\draw (\x,4) -- (\x,5.5);}
\end{scope}}

\newcommand{\brpos}[0]{
\draw (1,0) \br (0,1);
   \draw [white, line width=10](0,0) \br (1,1);
   \draw (0,0) \br (1,1);
}

\newcommand{\brrev}[0]{
\draw (0,0) \br (1,1);
\draw [white, line width=10](1,0) \br (0,1);
\draw (1,0) \br (0,1);}

\begin{document}
\begin{abstract}
We show that the ribbon zesting construction can produce {\it modular isotopes} -- different modular fusion categories with the same modular data. The result relies on the observation that the Reshetikhin-Turaev invariants of framed links associated to a ribbon fusion category satisfy a factorization property under zesting. This gives a new perspective on using topological invariants to classify topological order in light of modular data not being a complete invariant.

\end{abstract}
\maketitle

\section{Introduction}
The quantum invariants of just two framed links are the main tools used to classify modular fusion categories by rank. The $S$-matrix and twists that arise from coloring the Hopf link and unknot with $+1$ framing by simple objects are collectively referred to as the modular data. This data captures so much information about the underlying category that they are often treated as a proxy for the category itself.

However, while the modular data uniquely represents modular fusion categories at low rank, they are not a complete invariant in general. The first examples of modular fusion categories which are not determined by their modular data were discovered by Mignard and Schauenburg \cite{MS}. We call such categories {\it modular isotopes} for short.

\noindent \textbf{Definition.}
A set of inequivalent modular fusion categories are called \emph{modular isotopes} if they have the same modular data up to relabeling.

The name for these categories takes inspiration from their connection to condensed matter physics, where they serve as algebraic theories for quasiparticles called anyons and provide a mathematical approach to building a ``periodic table" for (2+1)D bosonic topological phases of matter. The modular data directly encodes several important characteristics of these quantum systems, namely the number of anyon types and their mutual and self-statistics. Even more information can be derived from the modular data, like the fusion/splitting rules, chiral central charge, etc.

While the Mignard-Schaunburg modular isotopes belong to a well-studied family of categories -- namely the twisted Drinfeld doubles of finite groups -- the smallest examples occur at rank 49, making them challenging to study directly without the aid of a computer. Following their discovery, several additional invariants beyond the modular data capable of distinguishing the Mignard-Schauenburg modular isotopes were identified. Some of these invariants are motivated as being in some sense minimal generalizations of the $S$-matrix \cite{BDGRTW, WW19, KMS} that are more discerning when used in place of the $S$-matrix together with the twists. However, an empirical analysis of how small knots and links perform as invariants for the smallest Mignard-Schauenburg modular isotopes identified many invariants with these properties \cite{DT}, albeit without any clear organizing principles.

Our main result is a new construction of Mignard and Schauenburg's modular isotopes that provides a framework for relating their topological invariants, by which we mean their quantum invariants of framed links and closed, oriented 3-manifolds. \\ \\
\noindent 
\textbf{Theorem \ref{thm:zestingmodularisotopes}.}
The Mignard-Schauenburg modular isotopes are related by cyclic Tannakian zesting with respect to their universal grading group $\mathds{Z}/p\mathds{Z}$.\\

Our proof relies on the observation that the Reshetikhin-Turaev invariants of framed links associated to a ribbon fusion category transform in a controlled way under the zesting construction. We show that for a given ribbon fusion category $\CC$, the framed link invariants associated to its zestings differ from those of $\CC$ by a scalar that can be computed directly from the input data $(\lambda, \nu, t,f)$ to the zesting construction.\\\\
\noindent \textbf{Theorem \ref{thm:zestinglinks}.}
Let $\CC$ be an $A$-graded ribbon fusion category, let $\mathcal{C}^{\lambda,\nu,t,f}$ be a ribbon zesting of $\CC$, and let $L$ be a framed $n$-component link whose $j$-th component is colored by $X_{i_j} \in \Irr(\CC_{i_j})=\Irr(\mathcal{C}^{\lambda,\nu,t,f}_{i_j})$. Then
\begin{align}
    L^{\lambda,\nu,t,f}_{X_{i_1}, X_{i_2}, \ldots, X_{i_n}}&=J(\lambda, \nu, t,f, i_1, i_2, \ldots, i_n)\cdot L_{X_{i_1}, X_{i_2}, \ldots, X_{i_n}}
\end{align}
where $J$ is a scalar-valued function of the ribbon zesting data $(\lambda, \nu, t,f)$ and the gradings $i_j$.\\

In other words, the Reshetikhin-Turaev invariant of a framed link $L$ factorizes under zesting. Moreover, our proof of Theorem \ref{thm:zestinglinks} shows it is straightforward to compute the factor $J$ for a given link. Such factorizations for the modular data -- and crucially, the Borromean rings invariant of the Mignard-Schauenburg modular isotopes established in \cite{KMS} -- are essential in proving Theorem \ref{thm:zestingmodularisotopes}.

Zesting then explains precisely the relationship between the topological invariants of the Mignard-Schauenburg modular isotopes.

The paper is structured as follows. In Section \ref{sec:preliminaries} we review the string diagram formalism for strict ribbon fusion categories and the ribbon zesting construction. In Section \ref{sec:zestingproperties} we define a product operation on ribbon zesting data that endows it with the structure of a monoid, which is helpful for later arguments. Section \ref{sec:zestinglinks} examines how topological invariants associated to ribbon and modular fusion categories transform under zesting and proves Theorem \ref{thm:zestinglinks}. The final Section \ref{sec:zestingmodularisotopes} begins with a review of twisted Drinfeld doubles of finite groups and the necessary details about Mignard-Schauenburg modular isotopes in particular. We then show that the Mignard-Schauenburg modular isotopes are related by zesting, proving Theorem \ref{thm:zestingmodularisotopes}.  Section \ref{sec:zestingmodularisotopes} ends with a brief interpretation of our results and places them in dialogue with the problem of classifying modular fusion categories. Several open questions are also presented.

\subsection*{Acknowledgements}

We thank Noah Snyder, Eric Samperton, Eric Rowell, Calvin McPhail-Snyder, and Constantin Teleman for helpful discussions. We also thank C\'{e}sar Galindo for pointing out a counterexample (see Example \ref{examp:cesar}) to a conjecture stated in an earlier version of this article. This work began with the Research Experiences for Undergraduates (REU) program at Indiana University supported by NSF grant \#1757857. The first author was supported by NSF Mathematical Sciences Postdoctoral Fellowship \#2002222. The third author was partially supported by NSF grants DMS-1802503 and DMS-1917319.

\section{Preliminaries} 
\label{sec:preliminaries}

In Section \ref{sec:ribboncats} we give a brief overview of the structure of a ribbon fusion category and its associated graphical calculus to set notation and terminology. We omit those details that aren't strictly necessary to understand the results of later sections and refer the reader to \cite{BK, EGNO} for formal definitions. Section \ref{sec:linkinvariants} recalls how the Reshetikhin-Turaev invariants of framed links arise from a ribbon fusion category. Lastly, Section \ref{sec:zesting} gives a concise introduction to the zesting construction on graded ribbon fusion categories with special attention to the case where the grading group is cyclic. 

\subsection{Ribbon categories and ribbon diagrams}
\label{sec:ribboncats}

It will be helpful to think of a ribbon fusion category in terms of layers of structure on a finite $\mathds{k}$-linear semisimple abelian category $\CC$. In general $\mathds{k}$ can be any algebraically closed field with characteristic zero, although in our applications we are interested in $\mathds{k}=\mathds{C}$. We write objects in the category as $X\in \CC$ and morphisms as $f \in \Hom(X,Y)$. 

\subsubsection{Fusion categories}
\label{sec:gradedfuscats}

The first layer of structure on $\CC$ is that of a coherent tensor product. This consists of a tensor product bifunctor $\otimes: \CC \times \CC \to \CC$ and a tensor unit object $\mathds{1}$ with natural isomorphisms $\alpha_{X,Y,Z}: (X \otimes Y) \otimes Z \to X \otimes (Y \otimes Z)$ for all $X,Y,Z \in \CC$ called associators satisfying the pentagon axioms \cite[Definition 2.1.1]{EGNO} and natural isomorphisms $l_X: \mathds{1} \otimes X \to X$, $r_X: X \otimes \mathds{1} \to X$ called left and right unitors satisfying the triangle axioms \cite[Definition 2.2.8]{EGNO}. In particular the tensor unit $\mathds{1}$ is a {\it simple} object, i.e.~ $\Hom(\mathds{1},\mathds{1}) \cong \mathds{C}$. We will write $\Irr(\CC)$ to refer to fixed set of representatives of isomorphism classes of simple objects. We write the {\it fusion rules} that dictate how the tensor product of representative simple objects decompose into simples as $X \otimes Y = \bigoplus_{Z \in \Irr(C)} N_{X,Y}^Z Z$. 

In the string diagram formalism we will draw morphisms $f \in \Hom(X,Y)$ from the top down, so that composition of morphisms is given by vertical stacking down the page. The tensor product of objects and morphisms is depicted by horizontal juxtaposition. To illustrate these conventions we show the picture that expresses compatibility of composition and tensor product, where the dashed lines indicate the two different ways of parsing the operations.

\begin{center}
\begin{tikzpicture}[line width=1, scale=.75, baseline=.125cm]
\onebox{}{X}{$f$}
\begin{scope}[yshift=-1.5cm]
\onebox{Z}{}{$g$}
\end{scope}
\begin{scope}[xshift=1.5cm]
\onebox{}{Y}{$h$}
\begin{scope}[yshift=-1.5cm]
\onebox{W}{}{$k$}
\end{scope}
\end{scope}
\draw[red,dashed] (.75,2) -- (.75,-1.5);
\end{tikzpicture}
\quad = \quad 
\begin{tikzpicture}[line width=1, scale=.75, baseline=.125cm]
\onebox{}{X}{$f$}
\begin{scope}[yshift=-1.5cm]
\onebox{Z}{}{$g$}
\end{scope}
\begin{scope}[xshift=1.5cm]
\onebox{}{Y}{$h$}
\begin{scope}[yshift=-1.5cm]
\onebox{W}{}{$k$}
\end{scope}
\end{scope}
\draw[red,dashed] (-.5,.25) -- (2,.25);
\end{tikzpicture}
\end{center}

In general the associators are depicted as a box on three strands, unless the category is {\it strict} and $\alpha_{X,Y,Z}=\id_{X\otimes Y\otimes Z}$, $l_X = \id_X$, and $r_X = \id_X$ for all $X,Y,Z \in \CC$, in which case strands can be freely grouped.

\begin{center}
\begin{tabular}{ccc}
\begin{tikzpicture}[line width=1,baseline=1cm]
\draw (0,0) node[below] {$X$} --(0,2) node[above] {$X$};
\draw (1,0) node[below] {$Y$}--(1,2) node[above] {$Y$};
\draw (2,0) node[below] {$Z$}--(2,2) node[above] {$Z$};
\draw[fill=white] (-.5,.75) rectangle node {$\alpha_{X,Y,Z}$} (2.5,1.25);
\end{tikzpicture}  & \qquad vs.\qquad &  \begin{tikzpicture}[line width=1, baseline=1cm]
\draw (0,0)node[below] {$X$}--(0,2) node[above] {$X$};
\draw (1,0) node[below] {$Y$} -- (1,2) node[above] {$Y$};
\draw (2,0)node[below] {$Z$} --(2,2)node[above] {$Z$};
\end{tikzpicture}\\
Generic associator & & Strict associator (identity morphism)\\
\end{tabular}
\end{center}

Every fusion category is equivalent to a strict fusion category \cite[Remark 2.8.7]{EGNO}. Since the graphical calculus on string diagrams in a strict fusion category has particularly nice topological properties that we will want to employ, from now on we draw string diagrams in a strict fusion category $\CC$ unless indicated otherwise. 

Fusion categories have the property of being {\it rigid}, guaranteeing (left) dual objects $X^*$ for all $X\in \CC$ and morphisms $X^* \otimes X \to \mathds{1}$ and $\mathds{1} \to X \otimes X^*$ called evaluation and coevaluation compatible with associativity. In pictures, the evaluation and coevaluation morphisms are drawn as cups and caps that satisfy a ``zig-zag axiom" and can be used to draw dual morphisms $f^* \in \Hom(Y^*, X^*)$ of $f\in \Hom(X,Y)$.
\begin{center}
\begin{tabular}{ccc}

 \begin{tikzpicture}[line width=1,xscale=.5,yscale=.5, baseline=.5cm,looseness=1.5]
       \draw[looseness=1.5] (0,0)--(0,2) to [out=90,in=90] (-1,2)--(-1,-.5) node[below] {$X$};
    \draw[looseness=1.5] (0,0) to [out=-90,in=-90] (1,0)--(1,2.5) node[above] {$X$};
   \end{tikzpicture} \quad = \quad \begin{tikzpicture}[line width=1,scale=.5, baseline=.5cm]
   \draw (0,-.5) node[below] {$X$}--(0,2.5) node[above] {$X$};
   \end{tikzpicture}
\quad
and
 \quad
   \begin{tikzpicture}[line width=1,xscale=-.5,yscale=.5, baseline=.5cm,looseness=1.5]
       \draw[looseness=1.5] (0,0)--(0,2) to [out=90,in=90] (-1,2)--(-1,-.5) node[below] {$X^*$};
    \draw[looseness=1.5] (0,0) to [out=-90,in=-90] (1,0)--(1,2.5) node[above] {$X^*$};
   \end{tikzpicture} \quad = \quad \begin{tikzpicture}[line width=1,scale=.5, baseline=.5cm]
   \draw (0,-.5) node[below] {$X^*$}--(0,2.5) node[above] {$X^*$};
   \end{tikzpicture} 
   & \qquad \qquad & 
   \begin{tikzpicture}[line width=1,scale=.75,baseline=.75cm]
\onebox{X^*}{Y^*}{f$^*$}
\end{tikzpicture}
\quad := \quad 
   \begin{tikzpicture}[line width=1,xscale=-.5,yscale=.5, baseline=.5cm,looseness=1.5]
       \draw[looseness=1.5] (0,0)--(0,2) to [out=90,in=90] (-1,2)--(-1,-.5) node[below] {$X^*$};
    \draw[looseness=1.5] (0,0) to [out=-90,in=-90] (1,0)--(1,2.5) node[above] {$Y^*$};
          \onebox{}{}{f}
   \end{tikzpicture}
 \\\\
  Rigidity axioms & & Dual morphisms 
   \end{tabular}
   \end{center}
The definition of rigidity also involves analogous right dual objects and associated evaluation and coevaluation morphisms \cite[Definition 2.10.11]{EGNO}. However, later we will recall that ribbon fusion categories are in particular {\it spherical}, and so for our purposes we can get away without explicitly referring to right duals.\\

Invariants of fusion categories that we will encounter include the rank $|\Irr(\CC)|$, fusion rules, and Frobenius-Perron dimensions. The Frobenius-Perron dimensions of simple objects $\FPdim(X)$ are the Frobenius-Perron eigenvalues of the {\it fusion matrices} $N_X$ for $X \in \Irr(\CC)$ where $(N_X)_{Y,Z} = N_{X,Y}^Z$. The Frobenius-Perron dimension of the category is $\FPdim(\CC)=\sum_{X \in \Irr(\CC)}\FPdim(X)^2$. 

The simplest examples of fusion categories are {\it pointed}, where every simple object $X \in \CC$ is invertible, i.e. $X \otimes X^* \cong \mathds{1}$.  Every pointed fusion category is equivalent to $\Vect_G^\omega$, the category of $G$-graded finite dimensional $\mathds{k}$-vector spaces for some group $G$ with associativity constraint determined by a $3$-cocycle $\omega$ representing a class in  $H^3(G, \mathds{k}^\times)$. The {\it pointed subcategory} $\CC_{pt}$ of a fusion category $\CC$ is the subcategory generated by invertible objects. The isomorphism classes of invertible objects always form a group under fusion, which we denote by $\Inv(\CC)$. 

A useful notion for working with fusion categories is to organize the fusion of objects using a finite group. A {\it $G$-grading} of a fusion category $\CC$ by a group $G$ is a decomposition $\CC = \bigoplus_{g \in G} \CC_g$, where $\CC_g$ are full abelian subcategories of $\CC$ such that $X_g \otimes Y_h \in \CC_{gh}$ for all $X_g \in \CC_g$, $Y_h \in \CC_h$. A grading is called faithful if $\CC_g \neq 0$ for all $g\in G$. Every fusion category has a {\it universal grading group}, i.e. a largest group by which it is faithfully graded. Every group which gives a faithful grading on the category is a quotient of the universal grading group \cite[Corollary 3.7]{GN}. The trivial component of the universal grading of $\CC$ coincides with the {\it adjoint subcategory} $\CC_{ad}$, the fusion subcategory of $\CC$ generated by $X \otimes X^*$ for simple $X$ in $\CC$.  

\subsubsection{Braided fusion categories}
\label{sec:braidedfuscats}
The next layer of structure is a {\it braiding} consisting of natural isomorphisms $c_{X,Y}: X\otimes Y \to Y \otimes X$ for all $X,Y \in \CC$ satisfying the hexagon axioms \cite[Definition 8.1.1]{EGNO}. We take the following convention for drawing braid isomorphisms and their inverses.

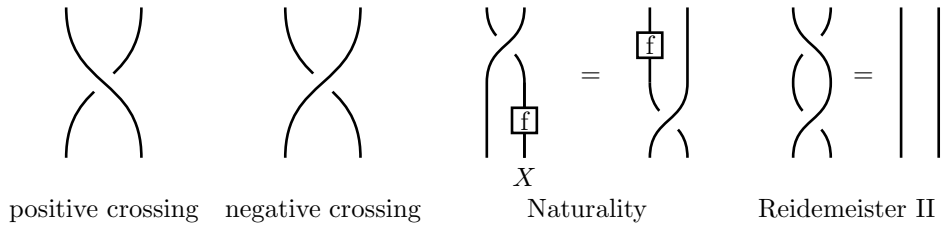
\begin{center}
\begin{tabular}{cc}
\begin{tikzpicture}[line width=1,scale=.5, looseness=1.25, baseline=1cm]
  \draw (0,0) node[below] {$Y$} \br (1,2)node[above] {$Y$};
 \draw[white, line width=10] (1,0) \br (0,2);
 \draw (1,0)node[below] {$X$} \br (0,2)node[above] {$X$};
\end{tikzpicture}
 \qquad &  \qquad 
\begin{tikzpicture}[line width=1,scale=.5, looseness=1.25, baseline=1cm]
 \draw (1,0) node[below] {$X$} \br (0,2) node[above] {$X$};
    \draw[white, line width=10] (0,0) \br (1,2);
    \draw (0,0) node[below] {$Y$} \br (1,2) node[above] {$Y$};
    \end{tikzpicture} 
\\
\text{Braiding} & \qquad \text{Inverse braiding}   \\
\end{tabular}
\end{center}
The axioms of a braided fusion category encode several important relations on string diagrams, like naturality and the Reidemeister II and III moves.

\begin{center}
\begin{tabular}{ccc}
\quad \begin{tikzpicture}[scale=.5, line width=1,looseness=1.25,baseline=1cm]
\draw (1,2) \br (0,4) node[above] {};
\draw[line width=10, white] (0,0)--(0,2) \br (1,4);
\draw (0,0) node[below] {} -- (0,2) \br (1,4) node[above] {};
\begin{scope}[xshift=1cm]
\onebox{}{}{$\scriptstyle g$}
\end{scope}
\onebox{}{}{$\scriptstyle f$}
\end{tikzpicture}
\quad
=
\quad
\begin{tikzpicture}[scale=.5, line width=1,looseness=1.25,baseline=1cm]
\draw (1,0) node[below] {} \br (0,2) node[above] {};
\draw [white, line width=10](0,0) \br (1,2) --(1,4);
\draw (0,0) node[below] {}\br(1,2)  -- (1,4) node[above] {};
\begin{scope}[yshift=2cm]
\onebox{}{}{$\scriptstyle g$}
\end{scope}
\begin{scope}[yshift=2cm,xshift=1cm]
\onebox{}{}{$\scriptstyle f$}
\end{scope}
\end{tikzpicture} \qquad &\qquad 
 \begin{tikzpicture}[scale=.5, line width=1,looseness=1.25,baseline=1cm]
    \draw (1,0) \br (0,2) \br (1,4);
    \draw[white, line width=10] (0,0) \br (1,2) \br (0,4);
    \draw (0,0) \br (1,2) \br (0,4);
    \end{tikzpicture}
    =\quad 
    \begin{tikzpicture}[scale=.5, line width=1,looseness=1.25,baseline=1cm]
    \draw (0,0) --(0,4);
    \draw (1,0) --(1,4);
    \end{tikzpicture} \qquad & \qquad
     \begin{tikzpicture}[xscale=.5,yscale=.34, line width=1,looseness=1.25,baseline=1cm]
  \draw (1,0) \br (0,2)--(0,4);
  \draw (2,0)--(2,2) \br (1,4) \br (0,6);
   \draw[white, line width=10] (0,0) \br (1,2) \br (2,4);
  \draw (0,0) \br (1,2) \br (2,4)--(2,6);
  \draw[white, line width=10] (0,4) \br (1,6);
  \draw (0,4) \br (1,6);
  \end{tikzpicture}
  \quad =\quad \begin{tikzpicture}[xscale=.5,yscale=.34, line width=1,looseness=1.25,baseline=1cm]
\draw (2,0) node[below] {} \br (1,2) \br (0,4)--(0,6) node[above] {};
\draw (2,4) \br (1,6);
\draw[line width=10, white] (0,0) -- (0,2) \br (1,4) \br (2,6);
\draw[line width=10, white] (1,0)\br (2,2) -- (2,4);
\draw (0,0) node[below] {} --(0,2)  \br (1,4) \br (2,6) node[above] {};
\draw (1,0) node[below] {}  \br (2,2)-- (2,4) node[above] {};
 \end{tikzpicture}\\
 Naturality & Reidemeister II & Reidemeister III
 \end{tabular}
\end{center}
  
If $c_{Y,X}=c_{X,Y}^{-1}$ for all $X, Y \in \CC$ then the category is called {\it symmetric}. Every symmetric fusion category is either {\it Tannakian} or {\it super-Tannakian}. Tannakian categories are of the form $\Rep(G)$ for $G$ a finite group and super-Tannakian categories are of the form $\Rep(G,z)$ where $z$ is a central element of order two. 

On the other end of the spectrum, if the only $X$ for which $c_{Y,X} \circ c_{X,Y}=\id_{X\otimes Y}$ for all $Y \in \CC$ is $X=\mathds{1}$, then the braiding and the category are called {\it nondegenerate}.

Not every fusion category $\CC$ admits a braiding, but one can always obtain a braided fusion category by taking its {\it Drinfeld center} $\mathcal Z(\CC)$ whose objects are pairs $(Z,\gamma)$ of objects $Z\in \CC$ and natural isomorphisms compatible with the monoidal structure of $\CC$ called {\it half-braidings} $ \gamma_X: X \otimes Z \to Z \otimes X $ for all $X \in \CC$ \cite[Definition 7.13.1]{EGNO}.

Note that when $\CC$ is braided, all gradings are by abelian groups $A$. We will mostly use additive notation for abelian groups, although there are a few places where it will be convenient to deviate.

\subsubsection{Ribbon fusion categories}

A braided fusion category $\CC$ is said to have a twist if there is a natural isomorphism $\theta$ of the identity functor $\id_\CC$ compatible with the braiding via $\theta_{X \otimes Y}=(\theta_X \otimes \theta_Y) \circ c_{Y,X} \circ c_{X,Y}$. If in addition $(\theta_X)^*=\theta_{X^*}$, then 
the twist is called a ribbon and $\CC$ is called a ribbon fusion category. We draw the twist isomorphisms and their inverses as follows.
\begin{center}
\begin{tikzpicture}[scale=.5, line width=1,looseness=1,baseline=.5cm]
\kink{X}{X}
\end{tikzpicture}
= 
$\theta_X$
\begin{tikzpicture}[scale=.5, line width=1,looseness=1,baseline=.5cm]
\draw(0,0) node[below] {$X$}--(0,2)node[above] {$X$};
\end{tikzpicture}
\qquad \qquad
\begin{tikzpicture}[yscale=-1, scale=.5, line width=1,looseness=1,baseline=-.5cm]
\kink{}{}
\draw (0,2) node[below] {$X$};
\draw (0,0) node[above] {$X$};
\end{tikzpicture}
= 
$\theta_X^{-1}$
\begin{tikzpicture}[scale=.5, line width=1,looseness=1,baseline=.5cm]
\draw(0,0) node[below] {$X$}--(0,2)node[above] {$X$};
\end{tikzpicture} 
\end{center}

The diagrams of the twists in a ribbon fusion category satisfy the framed Reidemeister I move. 
\begin{center}
    \begin{tabular}{c}
    $\begin{tikzpicture}[scale=.5, line width=1,looseness=1,baseline=.75cm]
\kink{X}{}
    \begin{scope}[yshift=3.5cm,yscale=-1]
    \kink{}{}
    \end{scope}
    \draw(0,3.5) node[above] {$X$};
    \end{tikzpicture} \quad = 
    \begin{tikzpicture}[scale=.5, line width=1,looseness=1,baseline=.75cm]
\draw(0,0) node[below] {$X$}--(0,3.5)node[above] {$X$};
    \end{tikzpicture}$
    \\
    (Framed) Reideimeister I
    \end{tabular}
\end{center}

The structure of ribbon twists on a braided fusion category is equivalent to a {\it spherical pivotal} structure with which it is useful to work simultaneously. In particular this structure defines a trace which gives an assignment of a number $\Tr(f) \in \Hom(\mathds{1}, \mathds{1}) \cong \mathbb{C}$ to every $f \in \Hom(X,X)$

\begin{align*}
\Tr(f) & =
    \begin{tikzpicture}[line width=1,scale=.5,baseline=.45cm]
    \onebox{}{}{f} 
    \draw[looseness=2] (0,2) to [out=90,in=90] (1,2) --(1,0) to [out=-90,in=-90] (0,0);
    \end{tikzpicture}
    =
      \begin{tikzpicture}[line width=1,scale=.5,baseline=.5cm]
    \onebox{}{}{f} 
    \draw[looseness=2] (0,2) to [out=90,in=90] (-1,2) --(-1,0) to [out=-90,in=-90] (0,0);
    \end{tikzpicture}\quad.
\end{align*}

Two ribbon fusion categories $\CC$ and $\DD$ are ribbon equivalent if there is a braided monoidal functor $(F,J)$ where $F:\CC \to \DD$ is an equivalence, $J_{X,Y}: F(X) \otimes F(Y) \to F(X\otimes Y)$ are natural isomorphisms compatible with associators and braiding via the commutative diagrams in \cite[Definition 2.4.1 and 8.1.7]{EGNO}, and $F(\theta)=\theta_F$. We will call a ribbon equivalence a {\it symmetry} of the category. Every symmetry of $\CC$ induces a (potentially trivial) permutation on $\Irr(\CC)$, which we refer to as a {\it relabeling}.

A nondegenerate ribbon fusion category is also known as a {\it modular} (fusion) category. 

\subsection{Reshetikhin-Turaev invariants of framed links and modular data}
\label{sec:linkinvariants}

Putting all of the pieces of the graphical calculus together (particularly isotopy and the framed Reidemeister moves), one has a recipe for assigning numerical invariants to framed links colored by simple objects of $\CC$: First represent the link $L$ as the Markov closure $\hat{b}$ of an $n$-strand braid $b \in B_n$ and color the strands of $b$ by objects $X_1, X_2, \ldots, X_n \in \CC$ in a manner that respects the labeling of the components of $L$. Then the invariant associated to $L$ is the trace of the morphism corresponding to the braid diagram. When working with an arbitrary link we will abuse notation and conflate the link $L$ with the value of its Reshetikhin-Turaev invariant and write
\begin{align*}
L_{X_1, X_2, \ldots, X_n} & =  \Tr(b)  = 
\begin{tikzpicture}[line width=1,scale=.5,baseline=.45cm]
\draw[looseness=2] (2,1.75) to [out=90,in=90] (4.5,1.75) --(4.5,.25) to [out=-90,in=-90] (2,.25)--(2,1.75);
\draw[looseness=2] (3,1.75) to [out=90,in=90] (4,1.75) --(4,0.25) to [out=-90,in=-90] (3,0.25)--(3,1.75);
\draw[fill=white] (1, .5) rectangle node {$b$} (3.5,1.5);
\draw (2.5,1.75) node {\scriptsize $\cdots$};
\draw (2.5,.25) node {\scriptsize $\cdots$};
\end{tikzpicture}\quad .
\end{align*}
Many algebraic invariants of ribbon fusion categories can be realized as such topological invariants. For example, the {\it categorical} or {\it quantum dimensions} 
\begin{align*} d_X:= \Tr(\id_X) =  X \, \,\,\begin{tikzpicture}[line width=1,scale=.3,baseline=-.125cm]
    \draw (0,0)  circle (1);
    \draw[->] (-1,.1)--(-1,-.1);
    \end{tikzpicture} &\quad  \end{align*} 
of simple objects $X \in \Irr(\CC)$ and {\it global quantum dimension} $D$ where $D^2=\sum_{X \in \Irr(\CC)} d_{X}^2$. (For {\it unitary} ribbon fusion categories, the quantum dimensions $d_X$ are guaranteed to be positive and coincide with the Frobenius-Perron dimensions $\FPdim(X)$ of the underlying fusion category. The global quantum dimension $D$ is also positive in the unitary case and is related to the Frobenius-Perron dimension of the category by $D^2=\FPdim(\CC)$.)

The {\it twists} $\theta_X$ of simples are realized as 
\begin{align*} \theta_X = \frac{1}{d_X}\Tr(c_{X,X}) = \frac{1}{d_X} \begin{tikzpicture}[baseline=0, line width=1,scale=.0675, shift={(0,-4.8)}]
\draw (2,2) to [out=0, in=180] (8,8);
\draw[white, line width=7](2,8) to [out=0, in=180] (8,2);
\draw(2,8) to [out=0, in=180] (8,2);
\draw [looseness=1](2,8) to [out=180, in=180] (2,2);
\draw [looseness=1](8,2) to [out=0, in=0] (8,8);
\draw (0,4) node[left] {\small $X$};
\draw[<-] (.25,5.1)--(.25,4.9);
\end{tikzpicture}.
\end{align*} 

Important information about the braiding is measured by the {\it $S$-matrix} associated to the Hopf link
\begin{align}
\label{smatrix}
    S_{X,Y}=\Tr\left(c_{Y^*,X} \circ c_{X,Y^*}\right) = \begin{tikzpicture}[baseline=0, line width=1,xscale=-1,scale=.2, shift={(0,-4.8)}]
                \draw (2, 5+1.25) arc (90:360:1.25cm);
                \draw (2, 5+1.25) arc (90:45:1.25cm);
                \draw (2+1.25, 5) arc (0:30:1.25cm);
                \draw (2-1.25,5) node[right] {\small $Y^*$};
                \draw (4-1.25, 5) arc (180:-90:1.25cm);
                \draw (4, 5-1.25) arc (270:225:1.25cm);
                \draw (4-1.25, 5) arc (180:210:1.25cm);
                \draw (7-1.25,5) node[left] {\small $X$};
                \draw[<-] (0.75,5.1)--(0.75,4.9);
                \draw[<-] (5.25,5.1)--(5.25,4.9);
\end{tikzpicture}.
\end{align}

The braiding on $\CC$ is nondegenerate (i.e.~$\CC$ is modular) if and only if the $S$-matrix is non-singular. For modular fusion categories the Reshetikhin-Turaev construction lifts to a (2+1)D TQFT and hence also yields 3-manifold invariants, but apart from a few brief comments (Remark \ref{rmk:zesting3manifoldinvariants1}) our focus here is on the framed link invariants. 

When $\CC$ is modular $S$ together with the diagonal {\it $T$-matrix} of twists 
\begin{align}
    \label{tmatrix} T_{X,Y} = \theta_{X}\delta_{X,Y}
\end{align} for $X,Y \in \Irr(\CC)$ satisfy the additional relations of the modular group $SL(2,\mathds{Z})$. For this reason when $\CC$ is modular the $S$-matrix and $T$-matrix together are called the {\it modular data}.

Many important pieces of data of a modular fusion category can be extracted from the modular data. The rank and global quantum dimension $D$ can be inferred from the dimensions and $S_{1,1}$, respectively. The categorical dimensions $d_X$ of simple objects can be read off of the first row (or column) of the $S$-matrix and the fusion rules can be found through the Verlinde formula \cite[Corollary 8.14.4]{EGNO}.  Other meaningful invariants that can be derived from modular data include the Frobenius-Schur indicators \cite{BDGRTW}, higher Gauss sums \cite{NSW}, and any link invariant which can arise from the closure of a 2-strand braid \cite{BDGRTW}. 

Due to their ability to capture so much information, the modular data often serve as a stand-in for the category itself and serve as the main tool for classifying modular fusion categories by rank \cite{RSW}. Moreover, in the first few decades of the study of modular fusion categories, all known examples were determined by their modular data. However, Mignard and Schauenburg eventually demonstrated a family of modular fusion categories for which the modular data are not a complete invariant \cite{MS}, which we now briefly describe.

Fix odd primes $p,q$ with $p>3$ and $p | q-1$. Then there are $p$ distinct Drinfeld centers $\ZZ(\Vect_{\mathds{Z}_q \rtimes \mathds{Z}_p}^{\omega})$ given by a choice of $3$-cocycle $\omega$ representing one of the $p$ classes in $H^3(G,\mathds{C}^{\times})$ but only $3$ distinct sets of modular data up to relabeling \cite{MS}. One set of modular data corresponds uniquely to the untwisted $\ZZ(\Vect_{\mathds{Z}_q \rtimes \mathds{Z}_p})$ and the remaining two sets of modular data each come from $(p-1)/2$ different $\ZZ(\Vect_{\mathds{Z}_q \rtimes \mathds{Z}_p}^{\omega})$. We call either of these two groups of $(p-1)/2$ categories the Mignard-Schauenburg modular isotopes, and refer to the whole family of $p$ categories $\ZZ(\Vect_{\mathds{Z}_q \rtimes \mathds{Z}_p}^{\omega})$ as the Mignard-Schauenburg categories.

Further properties of the categories $\ZZ(\Vect_{G}^{\omega})$ and the Mignard-Schauenburg modular isotopes in particular are discussed in more detail in Section \ref{sec:zestingmodularisotopes}.

\subsection{The ribbon zesting construction}
\label{sec:zesting}

The basic idea of zesting a ribbon fusion category is to modify the fusion rules by tensoring with invertible objects and then (if possible) to change the monoidal, braiding, and ribbon structure in a compatible way. At each stage of this construction -- associative zesting, braided zesting, and ribbon zesting -- whether the new category admits the given structure is governed by an algebraic obstruction theory \cite{DGPRZ}. We will now recall the data of a ribbon zesting along with the conditions they must satisfy. 

\begin{defn}[Adapted from \cite{DGPRZ}]
\label{def:zesting}
Let $\CC$ be an $A$-graded ribbon fusion category. The data of a {\it ribbon zesting} for $\CC$ is a tuple $(\lambda,\nu, t, j, f)$ comprised of
\end{defn}
\begin{enumerate}[label=(\alph*)]
    \item {\it associative zesting} data consisting of 
    \begin{itemize}
        \item a 2-cocycle $\lambda \in Z^2(A,\Inv(\CC_e))$, normalized so that $\lambda(e,a)=\lambda(a,e)=\mathds{1}$ for all $a \in A$,
        \item a 3-cochain $\nu \in C^3(A, \mathds{k}^{\times})$, normalized so that $\nu(e,b,c)=\nu(a,e,c)=\nu(a,b,e)$ for all $a,b,c \in A$ and identified with an isomorphism $$\nu_{a,b,c}: \lambda(a,b) \otimes \lambda(ab,c) \to \lambda(b,c) \otimes \lambda(a,bc)$$ which we draw as \begin{align*}
 \begin{tikzpicture}[line width=1,scale=1.25, baseline = 1.25cm]
  \foreach \x in {0,1}{
  \draw(\x,0)--(\x,2);
  }
  \draw[fill=white] (-.25,.75) rectangle node {$\nu_{a,b,c}$} (1.25,1.25) ;
   \draw (0,2) node[above] {$\lambda(a,b)$};
      \draw (1,2) node[above] {$\lambda(ab,c)$};    
  \draw (0,0) node[below] {$\lambda(b,c)$};
   \draw (1,0) node[below] {$\lambda(a,bc)$};
   \end{tikzpicture} \qquad & \qquad \text{ or simply } &  \begin{tikzpicture}[line width=1,scale=1.25, baseline = 1.25cm]
  \foreach \x in {0,1}{
  \draw(\x,0)--(\x,2);
  }
  \draw[fill=white] (-.25,.75) rectangle node {$a,b,c$} (1.25,1.25) ;
   \draw (0,2) node[above] {$(a,b)$};
      \draw (1,2) node[above] {$(ab,c)$};    
  \draw (0,0) node[below] {$(b,c)$};
   \draw (1,0) node[below] {$(a,bc)$};
   \end{tikzpicture}
   \end{align*}
        satisfying the {\it associative zesting condition} of Equation \ref{assoczestingcondition},
    \end{itemize}
    \item {\it braided zesting} data consisting of 
\begin{itemize}
    \item isomorphism $t(a,b): \lambda(a,b) \to \lambda(b,a)$ and
    \item a function $j:A \to \Aut_{\otimes}(\Id_{\CC})$ 
    drawn as 
    \begin{align}
\begin{tikzpicture}[line width=1,scale=1.25, baseline=1.25cm]
\onecirc{\lambda(b,a)}{\lambda(a,b)}{t(a,b)}
\end{tikzpicture}  \qquad \text{ or simply }\qquad \begin{tikzpicture}[line width=1,scale=1.25, baseline=1.25cm]
\onecirc{(b,a)}{(a,b)}{a,b}
\end{tikzpicture} \quad, \quad & \text{ and } \qquad \begin{tikzpicture}[line width=1,scale=1.25, baseline=1.25cm]
\onebox{X}{X}{$j_a$}
\end{tikzpicture}\quad.
\end{align} satisfying the {\it braided zesting conditions} of Equations \ref{braidzestingconditions},
\end{itemize}
    \item {\it twist zesting} data consisting of a function\footnote{While the twist zesting data $f(a)$ can also be interpreted and drawn as isomorphisms $f(a): X_a \to X_a$ it will simplify matters to treat it as a scalar in the equations that define a ribbon zesting.} $f:A \to \mathds{k}^{\times}$ satisfying the twist zesting condition of Equation \ref{twistzestingcondition} and the ribbon zesting condition of Equation \ref{ribbonzestingcondition}.
\end{enumerate}
Given ribbon zesting data $(\lambda,\nu,t,j,f)$ for a strict $A$-graded fusion category $\CC$, one can construct another $A$-graded fusion category $\CC^{\lambda,\nu, t,j,f}$ \cite{DGPRZ}. The associators, braiding, and twist of $\CC^{\lambda,\nu,t,j,f}$ can be constructed directly from that of $\CC$ and $(\lambda,\nu,t,j,f)$, although we will defer the details until we need them in Section \ref{sec:zestinglinkdiagrams}.

\subsubsection{Zesting conditions and properties}
\label{sec:notationassumptions}
Here we collect the conditions that must be satisfied by the zesting data in order for the zested category to achieve the desired layers of structures and properties. As indicated in the pictures in Definition \ref{def:zesting}, from time to time we will suppress the $\lambda$, $\nu$, and $t$ in the strand labels and coupons of diagrams when there is no risk of confusion.\\\\

\noindent \textbf{Associative Zesting Condition.}
\label{assoczestingcondition}
\begin{align}
 \scalebox{.9}{\begin{tikzpicture}[line width=1,scale=1,baseline=2.5cm]
  \foreach \x in {0,1,2}{
  \draw(\x,0)--(\x,4);
  }
  \draw[fill=white] (-.25,.75) rectangle node {$b,c,d$} (1.25,1.25) ;
  \draw[fill=white] (.75,1.75) rectangle node {$a,bc,d$} (2.25,2.25);
  \draw[fill=white] (-.25,2.75) rectangle node {$a,b,c$} (1.25,3.25) ;
  \draw (0,0) node[below] { $(c,d)$};
   \draw (1,0) node[below] {$(b,cd)$};
    \draw (2,0) node[below] {$(a,bcd)$};
    \draw (0,4) node[above] {$(a,b)$};
      \draw (1,4) node[above] {$(ab,c)$};
      \draw (2,4) node[above] {$(abc,d)$};
  \end{tikzpicture}}
  =
 \scalebox{.9}{\begin{tikzpicture}[line width=1,scale=1, baseline=2.5cm]
   \foreach \x in {0,1}{
  \draw(\x,0)--(\x,1);
  \draw(\x,2.5)--(\x,4);
  }
  \draw(2,0)--(2,4);
   \begin{scope}[xshift=0cm,yshift=2.5cm]
   \pic[braid /.cd] {braid=s_1^{-1}};
  \end{scope}
   \draw[fill=white] (.75,.75) rectangle node {$a,b,cd$} (2.25,1.25);
   \draw[fill=white] (.75,2.75) rectangle node {$ab,c,d$} (2.25,3.25);
  \draw (0,0) node[below] {$(c,d)$};
   \draw (1,0) node[below] {$(b,cd)$};
    \draw (2,0) node[below] {$(a,bcd)$};
    \draw (0,4) node[above] {$(a,b)$};
      \draw (1,4) node[above] {$(ab,c)$};   
      \draw (2,4) node[above] {$(abc,d)$};
  \end{tikzpicture}}
  \end{align}\\\\
  \noindent \textbf{Braided Zesting Condition.}
\begin{align}  
\label{braidzestingconditions}
\scalebox{.9}{\begin{tikzpicture}[line width=1,baseline=80,xscale=1]
\draw (0,2)--(0,4);
\twobox{(c,a)}{}{(ca,b)}{}{b,c,a}
 \begin{scope}[yshift=2cm]
 \onebox{}{}{$j_a$}
 \begin{scope}[xshift=1cm]
 \onecirc{}{}{a,bc}
 \end{scope}
 \end{scope}
  \begin{scope}[yshift=4cm]
  \twobox{}{(a,b)}{}{(ab,c)}{a,b,c}
 \end{scope}
\end{tikzpicture}}
=
\scalebox{.9}{\begin{tikzpicture}[line width=1,baseline=80,scale=1]
\draw (1,0) node[below] {$(ca,b)$}--(1,6) node[above] {$(ab,c)$};
\onecirc{(c,a)}{}{a,c}
 \begin{scope}[yshift=2cm]
\twobox{}{}{}{}{b,a,c}
 \end{scope}
  \begin{scope}[yshift=4cm]
  \onecirc{}{(a,b)}{a,b}
 \end{scope}
\end{tikzpicture}}, \quad \omega(a,b;c) \cdot
\scalebox{.9}{\begin{tikzpicture}[line width=1,baseline=80,scale=1]
\draw (0,2)--(0,4);
\twobox{(c,a)}{}{(ca,b)}{}{c,a,b}
 \draw (1.55,1.25) node[] {$^{-1}$};
 \begin{scope}[xshift=1cm,yshift=2cm]
 \onecirc{}{}{ab,c}
 \end{scope}
  \begin{scope}[yshift=4cm]
  \twobox{}{(b,c)}{}{(a,bc)}{a,b,c}
 \draw (1.55,1.25) node[] {$^{-1}$};
 \end{scope}
\end{tikzpicture}}
=
\scalebox{.9}{\begin{tikzpicture}[line width=1,baseline=80,scale=1]
\draw (1,0) node[below] {$(ca,b)$}--(1,6) node[above] {$(a,bc)$};
\onecirc{(c,a)}{}{a,c}
 \begin{scope}[yshift=2cm]
\twobox{}{}{}{}{a,c,b}
 \draw (1.55,1.25) node[] {$^{-1}$};
 \end{scope}
  \begin{scope}[yshift=4cm]
  \onecirc{}{(b,a)}{b,c}
 \end{scope}
\end{tikzpicture}}
\end{align}

In the second of Equations \ref{braidzestingconditions}, $\omega(a,b;c):=\omega(a,b)(c)$ where $\omega(a,b)=\chi_{\lambda(a,b)} \circ j_{ab} \circ j_a^{-1}
\circ j_b^{-1}$ and $\chi_{\lambda(a,b)}(X)\otimes \id_{\lambda(a,b)} = c_{\lambda(a,b),X}\circ c_{X,\lambda(a,b)}$.\\\\

\noindent \textbf{Twist Zesting Condition.}
\begin{align}
\label{twistzestingcondition}
 \scalebox{.9}{$f(ab)$\begin{tikzpicture}[line width=1,baseline=50]
 \draw (0,0) node [below] {$\parcen{a}$} \br (1,2);
  \draw (1,0) node [below] {$\parcen{b}$} \br (2,2);
  \draw (0,2.5) \br (2,4) node[above] {$(a,b)$};
  \draw[white, line width=10] (2,1) \br (0,2.5);
  \draw (2,1) \br (0,2.5);
  \draw[white, line width=10] (1,2) \br (0,4);
 \draw (1,2) \br (0,4)node [above] {$\parcen{a}$};
  \draw[white, line width=10] (2,2) \br (1,4);
 \draw (2,2) \br (1,4) node [above] {$\parcen{b}$};
 \begin{scope}[xshift=2cm,scale=.5]
 \kink{(a,b)}{}
 \end{scope}
\end{tikzpicture}}
= 
 \scalebox{.9}{$f(a)f(b)$ \begin{tikzpicture}[line width=1,,baseline=50]
 \draw (0,0)node[below] {$\parcen{a}$}--(0,2);
 \draw (1,0)node[below] {$\parcen{b}$}--(1,2);
 \begin{scope}[yshift=2cm]
 \onebox{}{\parcen{a}}{$j_b$}
 \end{scope}
 \begin{scope}[xshift=1cm,yshift=2cm]
 \onebox{}{\parcen{b}}{$j_a$}
 \end{scope}
  \begin{scope}[xshift=2cm]
  \onecirc{(a,b)}{}{b,a}
   \begin{scope}[yshift=2cm]
  \onecirc{}{(a,b)}{a,b}
  \end{scope}
  \end{scope}
  \end{tikzpicture}}
  \end{align}\\\\
  \noindent \textbf{Ribbon Zesting Condition.}
  \begin{align}
\label{ribbonzestingcondition}
 f(a) = f(-a)\chi_{\lambda(a,-a)}(X_a)\theta_{\lambda(a,-a)}
 \end{align}

It will also be important to track how the trace and modular data transform under zesting. 

\begin{prop}[Proposition 5.5, Theorem 5.7 of \cite{DGPRZ}]
\label{prop:zestedtraceandmodulardata}
Let $(\lambda,\nu,t,j,f)$ be zesting data for an $A$-graded ribbon fusion category $\CC$. Then the trace in $\CC^{\lambda,\nu,t,j,f}$ is given by

 \begin{align*}\Tr^{(\lambda, \nu, t, j, f)}(\rho)= \frac{f(a)}{\dim(\lambda(-a,a)) t(a,a)}\Tr(j_a^{-1}(X_a)\circ \rho) \quad , \quad \rho \in \Hom_{\CC^{\lambda,\nu, t,j,f}}(X_a,X_a)
 \end{align*}
 and the $S$- and $T$-matrices by 

\begin{align}
\label{eq:zestedmodulardata}
S^{\lambda,\nu,t,j,f}_{X_a,Y_b} &= \frac{\dim(\lambda(a-b,b-a))\dim(\lambda(a,-b))t(a,-b)t(-b,a)f(a-b)}{t(a-b,a-b)\dim(X_a)}\frac{j_a(Y_b)j_b(X_a)}{j_{ab}(X_a)j_{ab}(Y_b)j_{ab}(\lambda(a,b))} S_{X_a,Y_b}S_{X_a,\lambda(b,-b)}\\
T^{\lambda,\nu,t,j,f}_{X_a,Y_b} &= f(a)\theta_{X_a} \delta_{X_a,Y_b}
\end{align}
 \end{prop}
 \qed
 
\begin{rmk}
Nondegeneracy of the braiding is not preserved by zesting in general. However, the zestings of a modular fusion category with respect to its universal grading group will remain modular. 
\end{rmk}

When the ribbon zesting uses the universal grading group we may assume that the $j$ in the zesting data is trivial. Recall from \cite[Definition 4.6]{DGPRZ} that for a fixed associative zesting $(\lambda, \nu)$, two braided zestings $(\lambda, \nu, t,j)$ and $(\lambda, \nu, t', j')$ are called {\it similar} if the braiding isomorphisms on $\CC^{\lambda,\nu,t,f}$ and $\CC^{\lambda,\nu,t',f'}$ are equal.
\begin{prop}[Corollary 4.8 of \cite{DGPRZ}]
\label{prop:jtrivial}
Let $\CC$ be an $A$-graded braided fusion category. When $A$ is the universal grading group of $\CC$ every zesting $(\lambda, \nu, t, j)$ is similar to a braided zesting of the form $(\lambda, \nu, t, \id_{\mathcal{C}})$. \qed
\end{prop}

Since our eventual application of zesting in Section \ref{sec:zestingmodularisotopes} is with respect to a universal grading group, from now on we will always assume that $j$ is trivial. We will also be working with $\CC$ unitary, so that all categorical dimensions are positive. 

Assuming that $j$ is trivial gives an even simpler formula for the zested trace than the one in Proposition \ref{prop:zestedtraceandmodulardata}, namely
\begin{align}
\label{eq:zestedtracejtrivial}\operatorname{Tr}^{(\lambda, \nu, t, f)}(\rho)=\frac{ f(a)}{ t(a,a)}\operatorname{Tr}(\rho).
\end{align}

As a result, when computed a zested trace one may substitute the with the original trace after correcting with the necessary scalar factors of $\frac{f(a)}{t(a,a)}$. We will revisit this fact in Section \ref{sec:zestinglinkdiagrams}.

\subsubsection{Cyclic Tannakian ribbon zesting}
\label{sec:cyclictannakianzesting}

When the 2-cocycle $\lambda$ takes values in a cyclic Tannakian subcategory $\Rep(\mathds{Z}/ N\mathds{Z})$ of $\CC$, the data of a $\mathds{Z}/ N\mathds{Z}$-ribbon zesting takes a particularly simple form. Fix a generator $g$ of $\Rep(\mathds{Z}/ N\mathds{Z})$ inside $\CC$, $q$ an $N$-th root of unity, and $\zeta$ a $2N$-th root of unity. Then the ribbon zestings are parametrized by triples $(a,b,s)$ where $a,b  \in \mathds{Z}/ N\mathds{Z}$ with $a = 
2b \mod N$ where 
$s$ is a solution to $s^N=\zeta^{-
2b
}$.

For a fixed $(a,b,s)$ the zesting data $(\lambda_a,\nu_b, t_s,f_s)$ takes the following form, adapted here from Section 6.2.3 and in particular Propositions 6.3-6.4 of \cite{DGPRZ}. 
\begin{align}
   \lambda_a(i,j) &= \begin{cases} \mathds{1} & i +j < N \\
    g^a & i +j \ge N \end{cases} \\
\nu_b(i,j,k) &= \begin{cases} 1 & i +j < N \\
\label{eq:3cochain}
   \zeta^{k
   2b
   }  & i+j \ge N \end{cases}\\
   \label{eq:cycliczestingt}
   t_s(i,j) & = s^{-ij}\id_{\lambda_a(i,j)} \\
    \label{eq:cycliczestingf}
   f_s(i) &=s^{-i^2}.
\end{align}

The following observation will let us simplify some later calculations.
\begin{prop}
\label{prop:cyclictrace}
Let $(\lambda_a,\nu_b, t_s, f_s)$ be a cyclic Tannakian zesting. Then 
\begin{align*} \Tr^{\lambda,\nu,t,f}(\rho) = \Tr(\rho) & \quad \rho \in \Hom(X,X), X \in \CC.
\end{align*}
\end{prop}
\begin{proof}
Let $X\in\CC_i$ be a simple object. By Equation \ref{eq:zestedtracejtrivial}, 
\begin{align*}
    \Tr^{\lambda,\nu,t,f}(\rho) &= \frac{f(i)}{t(i,i)} \Tr(\rho) =\frac{s^{-i^2}}{s^{-i^2}} \Tr(\rho)  = \Tr(\rho).
\end{align*}
\end{proof}

Finally we recall how the modular data transform under a cyclic Tannakian ribbon zesting. 

\begin{prop}[Proposition 6.5 of \cite{DGPRZ}]
\label{prop:cyclictannakianmodulardata}
Let $\CC$ be a $\mathds{Z}/N\mathds{Z}$-graded ribbon fusion category with zesting data $(\lambda_a,\nu_b, t_s,f_s)$ as above. Then the modular data of $\CC^{\lambda_a,\nu_b, t_s,f_s}$ are related to the modular data of $\CC$ by 
\begin{align}\label{eq:cyclictannakianSmatrix}
S_{X_i,Y_j}^{\lambda_a,\nu_b,t_s,f_s}=s^{2ij}S_{X_i,Y_j}, &&  0\leq i ,j < N,  X_i\in \CC_i, Y_j\in \CC_j,
\end{align}
\begin{align}\label{eq:cyclictannakianTmatrix}
T_{X_i,X_i}^{\lambda_a,\nu_b,t_s,f_s}= s^{-i^2}T_{X_i,X_i}, &&  0\leq i < N, X_i\in \CC_i,
\end{align}
where $s$ is a choice of solution to the equation $s^N=\zeta^{-
2b
}$ for pairs $a,b \in \mathds{Z}/N\mathds{Z}$ with $a=
2b \mod N$.\qed
\end{prop}

\section{Properties of ribbon zesting}
\label{sec:zestingproperties}

In this section we exhibit some properties of ribbon zesting that will be used later to relate the Mignard-Schauenburg modular isotopes. First we define a product operation on zesting data and show that this defines the structure of a monoid on the choices of (braided or) ribbon zesting data for a fixed category. We also point out a natural notion of inverse for zesting data under this product, although we show it is only a true inverse after passing to suitable equivalence classes.

\subsection{A product on ribbon zesting data}
\label{sec:seqzesting}
Fix a grading group $A$ of a (strict) ribbon fusion category $\CC$ and let $(\lambda,\nu,t,f)$ and $(\lambda',\nu',t',f')$ be two choices of ribbon zesting data satisfying Definition \ref{def:zesting}. We will show that there is a well-defined product of ribbon zesting data which we denote by $(\lambda \cdot \lambda', \nu \cdot \nu', t \cdot t', f \cdot f')$. 

\subsubsection{Product of associative zesting data}

The map 
\begin{align*} \lambda \cdot \lambda': A \times A &\to \Inv(\CC_e) \\
(a,b) & \mapsto \lambda(a,b) \otimes \lambda'(a,b)
\end{align*}
automatically satisifies the normalized 2-cocycle condition in the fusion ring of $\CC$ and defines a new fusion rule  $X_a \overset{\lambda\lambda'}{\otimes} Y_b := X_a \otimes Y_b  \otimes \lambda(a,b)
\otimes \lambda'(a,b)$. This fusion rule can then be seen to categorify to a tensor product bifunctor on $\CC$ via the isomorphism
\[(\nu \cdot \nu')_{a,b,c}: \lambda(a,b) \otimes \lambda'(a,b) \otimes  \lambda (ab,c) \otimes  \lambda '(ab,c) \to  \lambda(b,c) \otimes \lambda'(b,c)\otimes  \lambda(a,bc) \otimes \lambda'(a,bc) \]
defined by the picture
\begin{align}
\label{def:seqassoc}
\begin{tikzpicture}[line width=1, scale=1.85,baseline=50]
\cocycm{a}{b}{c}{(\lambda\cdot\lambda')}{(\nu\cdot \nu')}
\end{tikzpicture} && :=
 &  &
\begin{tikzpicture}[line width=1, baseline=60, xscale=1.25, yscale=.75]
\draw (1,.5) \br (2,2);
\draw[white, line width=10] (2,.5) \br (1,2);
\draw (2,.5) \br (1,2);
\begin{scope}[yshift=3.5cm]
\draw (2,.5) \br (1,2);
\draw[white, line width=10] (1,.5) \br (2,2);
\draw (1,.5) \br (2,2);
\end{scope}
\begin{scope}[yshift=2cm]
\twobox{}{}{}{}{\nu_{a,b,c}}
\end{scope}
\begin{scope}[xshift=2cm, yshift=2cm]
\twobox{}{}{}{}{\nu'_{a,b,c}}
\end{scope}
\foreach \x in {0,3}
{\draw (\x, .5)--(\x,2);
\draw (\x,4) -- (\x,5.5);}
\draw (0,5.5) node[above] {$\lambda(a,b)$};
\draw (1,5.5) node[above] {$\lambda'(a,b)$};
\draw (2,5.5) node[above] {$\lambda(ab,c)$};
\draw (3,5.5) node[above] {$\lambda'(ab,c)$};
\draw (0,0.5) node[below] {$\lambda(b,c)$};
\draw (1,0.5) node[below] {$\lambda'(b,c)$};
\draw (2,0.5) node[below] {$\lambda(a,bc)$};
\draw (3,0.5) node[below] {$\lambda'(a,bc)$};
\end{tikzpicture}.
\end{align}

Note that the product isomorphism $\nu \cdot \nu'$ is only the tensor product $\nu \otimes \nu'$ up to some braidings which must be introduced in order to identify the source and target objects.

Checking that $\nu\cdot \nu'$ satisfies the associative zesting condition amounts to proving the following picture equation.
\begin{align}
\begin{tikzpicture}[line width=1,baseline=3.75cm]
\draw (0,2.75)--(0,5.25);
\draw (1,2.75)--(1,5.25);
\draw (4,5.25)--(4,7.75);
\draw (5,5.25)--(5,7.75);
\draw (4,0.25)-- (4,3);
\draw (5,0.25)-- (5,3);
\zestedseqcochain{b}{c}{d}
\begin{scope}[xshift=2cm,yshift=2.5cm]
\zestedseqcochain{a}{bc}{d}
\end{scope}
\begin{scope}[yshift=5cm]
\zestedseqcochain{a}{b}{c}
\end{scope}
\end{tikzpicture}
\qquad =
\quad
\begin{tikzpicture}[line width=1,baseline=3.75cm]
\begin{scope}[xshift=2cm]
\zestedseqcochain{a}{b}{cd}
\end{scope}
\begin{scope}[xshift=2cm,yshift=5cm]
\zestedseqcochain{ab}{c}{d}
\end{scope}
\draw (0,.25) \br (2,5.25);
\draw (1,.25) \br (3,5.25);
\draw[white, line width=10] (2,2.75) \br (0,7.75);
\draw[white, line width=10] (3,2.75) \br (1,7.75);
\draw  (2,2.75) \br (0,7.75);
\draw (3,2.75) \br (1,7.75);
\draw (4,2.75)--(4,5.25);
\draw (5,2.75)--(5,5.25);
\draw[red,dashed] (2.25,2) rectangle (4.5,6);
\end{tikzpicture}.
\end{align}
Observe that on the left hand side the diagram can be isotoped into two disjoint diagrams, one involving 3 $\nu$ boxes and the other involving 3 $\nu'$ boxes. After applying a Reidemeister III move to the region in the red rectangle on the right hand side, the same can be said of the diagram on the right. The two sides can then be identified by applying the associative zesting conditions Equation \ref{assoczestingcondition} satisfied by $\nu$ and $\nu'$.

\begin{rmk}
In \cite[Definition 3.2]{DGPRZ} a more general notion of associative zesting was defined for $G$-graded {\it fusion} i.e. not necessarily braided categories.
This same argument shows that the product of such associative zestings is again an associative zesting. One need only take into account that the gradings no longer live in an abelian group and that the braidings instead depict half braidings in the appropriate relative center.
\end{rmk}

\subsubsection{Product of braided zesting data}

Now consider the isomorphism $$(t\cdot t')(a,b): = t(a,b) \otimes t'(a,b): \lambda(a,b) \otimes \lambda'(a,b) \to \lambda(b,a) \otimes \lambda'(b,a)$$ depicted by

\begin{align}
\begin{tikzpicture}[line width=1,scale=1.25,baseline=1.25cm]
\onecirc{\lambda''(b,a)}{(\lambda\cdot\lambda')(a,b)}{\scriptstyle (t\cdot t')(a,b)}
\end{tikzpicture}
&\hspace{10pt}:=\hspace{10pt}
\begin{tikzpicture}[line width=1,scale=1.25,baseline=1.25cm]
\onecirc{\lambda(b,a)}{\lambda(a,b)}{t(a,b)}
\begin{scope}[xshift=1.5cm]
\onecirc{\lambda'(b,a)}{\lambda'(a,b)}{t'(a,b)}
\end{scope}
\end{tikzpicture}\quad.
\end{align}

The first of the two braided zesting conditions of Equations \ref{braidzestingconditions} that must be satisfied by $t\cdot t'$ manifests as the picture equation

\begin{align}
\begin{tikzpicture}[line width=1,scale=.75, baseline=105pt]
\begin{scope}[yscale=.75]
\seqassoc{b}{c}{a}
\end{scope}
\begin{scope}[yshift=4.125cm,xshift=2cm]
\onecirc{}{}{t}
\begin{scope}[xshift=1cm]
\onecirc{}{}{t'}
\end{scope}
\end{scope}
\begin{scope}[yshift=5.75cm, yscale=.75]
\seqassoc{a}{b}{c}
\end{scope}
\draw (0,4.125)--(0,6.125);
\draw (1,4.125)--(1,6.125);
\end{tikzpicture}
\quad =\quad &
\begin{tikzpicture}[line width=1,scale=.75, baseline=105pt]
\draw (0,0.375)--(0,3.125);
\draw (1,.375)--(1,3.125);
\draw (0,6.875)--(0,9.875);
\draw (1,6.875)--(1,9.875);
\begin{scope}[yshift=.375cm]
\onecirc{}{}{t}
\end{scope}
\begin{scope}[xshift=1cm, yshift=.375cm]
\onecirc{}{}{t'}
\end{scope}
\begin{scope}[yshift=2.75cm, yscale=.75]
\seqassoc{b}{a}{c}
\end{scope}
\begin{scope}[yshift=7.875cm]
\onecirc{}{}{t}
\begin{scope}[xshift=1cm]
\onecirc{}{}{t'}
\end{scope}
\end{scope}
\draw (2,0.375)--(2,3.125);
\draw (3,.375)--(3,3.125);
\draw (2,6.875)--(2,9.875);
\draw (3,6.875)--(3,9.875);
\end{tikzpicture} \quad , 
\end{align}
where we have omitted the arguments for the $t$ isomorphisms for readability. Applying the first braided zesting condition for each of $\nu$ and $\nu'$ gives the identification of the two diagrams up to isotopy. Verifying the second braided zesting condition is a nearly identical calculation after one observes that the factor $\omega$ satisfies $(\omega\cdot\omega')(a,b;c) = \omega(a,b;c)\omega'(a,b;c)$ due to the fact that $\chi_{\lambda(a,b)}$ is a tensor natural isomorphism \cite[Proposition 2.8(i)]{DGPRZ}.

Once again, we remark that this defines a product on braided (not necessarily ribbon) zesting data $(\lambda,\nu,t)$.

\subsubsection{Product of twist zesting data}

Finally, we check that the product of twist zestings is a twist zesting. Given twist zestings $f: A \to \mathds{k}^{\times}$ and $f':A \to \mathds{k}^{\times}$consider  the pointwise product
\begin{align}
  f\cdot f': A \to \mathds{k}^{\times}.
\end{align}

After applying the definitions of $t \cdot t'$ and $f\cdot f'$ one uses that $\chi$, $\theta$, and $j$ are all {\it tensor} natural isomorphisms to show 
that the twist zesting condition Equation \ref{twistzestingcondition} factors into the two twist zesting equations satisfied independently by $f$ and $f'$. A similar argument shows that if $f$ and $f'$ are ribbon, then so is $f''$. 

\begin{lem}
\label{lem:seqzesting}
The choices of zesting data for a fixed strict (braided or) ribbon fusion category form a monoid under the product defined above.  
\end{lem}
\begin{proof}
The discussion preceding this lemma shows that the product of zesting data is well-defined. By assumption, the left and right unitors in $\CC$ are trivial. It follows that $(\lambda, \nu, t,f)=(\mathds{1}, \id, \id, 1)$ is an identify for the product operation on zesting data. To see that the product is associative, we check that $(\nu \cdot \nu')\cdot \nu''= \nu \cdot (\nu' \cdot \nu'')$ by verifying that the left and right hand sides below are isotopic.

\begin{align*}
\begin{tikzpicture}[line width=1,scale=.75, baseline=0]
\draw (0,-2)node [below] {$\scriptstyle \lambda(b,c)\phantom{'}$}--(0,2) node[above] {$\scriptstyle \lambda(a,b)$};
\draw (6,-2)node [below] {$\scriptstyle \lambda''(a,bc)$}--(6,2) node[above] {$\scriptstyle \lambda''(ab,c)$};
\draw (2,.25) to [out=90,in=-90] (1,2) node[above] {$\scriptstyle \lambda'(a,b)$};
\draw (5,.25) to [out=90,in=-90] (2,2);
\draw (2.125,2) node[above] {$\scriptstyle \lambda''(a,b)$};
\draw (2,-.25) to [out=-90,in=90] (1,-2) node [below] {$\scriptstyle \lambda'(b,c)$};
\draw (5,-.25) to [out=-90,in=90] (2,-2);
\draw (2.125,-2)  node [below] {$\scriptstyle \lambda''(b,c)$};
\draw[white,line width=10] (1,.25) to [out=90,in=-90] (3,2);
\draw (1,.25) to [out=90,in=-90] (3,2);
\draw (3.25,2) node[above] {$\scriptstyle \lambda(ab,c)$};
\draw[white,line width=10] (1,-.25) to [out=-90,in=90] (3,-2);
\draw (1,-.25) to [out=-90,in=90] (3,-2);
\draw (3.375,-2) node[below] {$\scriptstyle \lambda(a,bc)\phantom{'}$};
\draw[white,line width=10] (3,.25) to [out=90,in=-90] (4,2);
\draw (4.5,2) node[above] {$\scriptstyle \lambda'(ab,c)$};
\draw (3,.25) to [out=90,in=-90] (4.25,2);
\draw[white,line width=10] (3,-.25) to [out=-90,in=90] (4.25,-2);
\draw (3,-.25) to [out=-90,in=90] (4,-2);
\draw (4.675,-2) node[below] {$\scriptstyle \lambda'(a,bc)$};
\draw[fill=white] (-.25,-.25) rectangle node {$\nu$} (1.25,.25);
\begin{scope}[xshift=2cm]
\draw[fill=white] (-.25,-.25) rectangle node {$\nu'$} (1.25,.25);
\end{scope}
\begin{scope}[xshift=5cm]
\draw[fill=white] (-.25,-.25) rectangle node {$\nu''$} (1.25,.25);
\end{scope}
\end{tikzpicture}
=
\begin{tikzpicture}[line width=1,scale=.75, baseline=0]
\draw (0,-2)node [below] {$\scriptstyle \lambda(b,c)\phantom{'}$}--(0,2) node[above] {$\scriptstyle \lambda(a,b)$};
\draw (6,-2)node [below] {$\scriptstyle \lambda''(a,bc)$}--(6,2) node[above] {$\scriptstyle \lambda''(ab,c)$};
\draw (3,.25) to [out=90,in=-90] (1,2) node[above] {$\scriptstyle \lambda'(a,b)$};
\draw (5,.25) to [out=90,in=-90] (2,2);
\draw (2.125,2) node[above] {$\scriptstyle \lambda''(a,b)$};
\draw (3,-.25) to [out=-90,in=90] (1,-2) node [below] {$\scriptstyle \lambda'(b,c)$};
\draw (5,-.25) to [out=-90,in=90] (2,-2);
\draw (2.25,-2) node [below] {$\scriptstyle \lambda''(b,c)$};
\draw[white,line width=10] (1,.25) to [out=90,in=-90] (3,2);
\draw (1,.25) to [out=90,in=-90] (3,2);
\draw (3.25,2) node[above] {$\scriptstyle \lambda(ab,c)$};
\draw[white,line width=10] (1,-.25) to [out=-90,in=90] (3,-2);
\draw (1,-.25) to [out=-90,in=90] (3,-2);
\draw (3.5,-2) node[below] {$\scriptstyle \lambda(a,bc)\phantom{'}$};
\draw[white,line width=10] (4,.25) to [out=90,in=-90] (5,2);
\draw (4,.25) to [out=90,in=-90] (5,2);
\draw (4.675,2) node[above] {$\scriptstyle \lambda'(ab,c)$};
\draw[white,line width=10] (4,-.25) to [out=-90,in=90] (5.,-2);
\draw (4,-.25) to [out=-90,in=90] (5,-2);
\draw (4.675,-2) node[below] {$\scriptstyle \lambda'(a,bc)$};
\draw[fill=white] (-.25,-.25) rectangle node {$\nu$} (1.25,.25);
\begin{scope}[xshift=3cm]
\draw[fill=white] (-.25,-.25) rectangle node {$\nu'$} (1.25,.25);
\end{scope}
\begin{scope}[xshift=5cm]
\draw[fill=white] (-.25,-.25) rectangle node {$\nu''$} (1.25,.25);
\end{scope}
\end{tikzpicture}
\end{align*}
The associativity of the braided and twist zesting data is similarly clear.
\end{proof}

We note that the product we have defined is not obviously commutative.

\subsection{Inverse zesting data}
Next we show that zestings can be ``inverted" up to equivalence. Let $(\lambda,\nu,t,f)$ be zesting data for $\CC$ and define $(\tilde{\lambda},\tilde{\nu},\tilde{t},\tilde{f})$ where $\tilde{\lambda}(a,b)=\lambda(a,b)^*$ and 

\begin{align}
    \tilde{\nu}_{a,b,c} &= \begin{tikzpicture}[line width=1,scale=.75, baseline=.5cm]
    \draw (0,0)--(0,2);
    \draw (1,0)--(1,2);
    \draw[fill=white] (-.25,.75) rectangle (1.25,1.25);
    \draw (.5,1) node {\small $a,b,c$};
    \draw (1.25,1.25) node[right] {$^{-1}$};
    \draw[looseness=1.5] (0,2) to [out=90,in=90] (3,2) -- (3,-1);
    \draw[looseness=1.5] (1,2) to [out=90,in=90] (2,2) -- (2,-1);
    \draw[looseness=1.5] (0,0) to [out=-90,in=-90] (-1,0) -- (-1,3);
    \draw[looseness=1.5] (1,0) to [out=-90,in=-90] (-2,0) -- (-2,3);
    \draw (3,-2.5) \br (2,-1);
     \draw[white, line width=10]  (2,-2.5) \br (3,-1);
     \draw  (2,-2.5) \br (3,-1);
    \draw (-1,3) \br (-2,4.5);
    \draw[white, line width=10] (-2,3) \br (-1,4.5);
      \draw (-2,3) \br (-1,4.5);
    \draw (1.5,-2.5) node[below] {$\lambda(b,c)^*$};
    \draw (3.5,-2.5) node[below] {$\lambda(a,bc)^*$};
    \draw (-2.5,4.5) node[above] {$\lambda(a,b)^*$};
    \draw (-.5,4.5) node[above] {$\lambda(ab,c)^*$};
    \end{tikzpicture} 
  &
    \tilde{t} &=  \begin{tikzpicture}[line width=1,scale=1,baseline=1cm]
    \onecirc{}{}{a,b}
    \draw (0.33,1.25) node[right]{$^{-1}$};
    \draw[looseness=1.5] (0,2) to [out=90,in=90] (1,2)--(1,-1) node[below] {$\lambda(b,a)^*$};
    \draw[looseness=1.5] (0,0) to [out=-90,in=-90] (-1,0)--(-1,3) node[above] {$\lambda(a,b)^*$};
    \end{tikzpicture}
   & 
    \tilde{f}&=f^{-1}.
\end{align}

It is straightforward to check that the choices of $\tilde{\lambda}$, $\tilde{\nu}$, $\tilde{t}$, and $\tilde{f}$ satisfy the conditions for zesting data in Equations \ref{assoczestingcondition} through \ref{ribbonzestingcondition}. Then by Lemma \ref{lem:seqzesting}, the composition of zestings  $\CC^{\lambda \cdot \tilde{\lambda}, \nu \cdot \tilde{\nu}, t \cdot \tilde{t}, f \cdot \tilde{f}}$ gives a well-defined ribbon fusion category, which we now check is ribbon-equivalent to $\CC$.

Recall from Section \ref{sec:ribboncats} that a ribbon equivalence is a braided monoidal equivalence $(F,J) $ which respects twists. To see that $\CC^{\lambda \cdot \tilde{\lambda}, \nu \cdot \tilde{\nu}, t \cdot \tilde{t}, f \cdot \tilde{f}} \simeq \CC$, let $F$ be the identity functor on $\CC^{\lambda \cdot \tilde{\lambda}, \nu \cdot \tilde{\nu}, t \cdot \tilde{t}, f \cdot \tilde{f}}$ and let $J_{X_a, Y_b}: X_a \otimes Y_b \otimes \lambda(a,b) \otimes \lambda(a,b)^* \to X_a \otimes Y_b$ be the isomorphism given by ``cupping off" the dual pair of invertible objects. 

\begin{align} F &= \Id_\CC, &  J_{X_a,Y_b} = 
\begin{tikzpicture}[line width=1, baseline=1cm]
\draw (0,0) node [below] {$X_a$}--(0,2) node [above] {$X_a$};
\draw (1.25,0) node [below] {$Y_b$}--(1.25,2) node [above] {$Y_b$};
\draw[looseness=2] (2.5,2) node [above] {$\lambda(a,b)$} to [out=-90,in=-90] (3.75,2) node [above] {$\lambda(a,b)^*$};
\end{tikzpicture}
\end{align}

The axiom that establishes $(\Id_{\CC},J)$ as a monoidal equivalence is satisfied if the following picture equation holds.
\begin{align}
\begin{tikzpicture}[line width=1,scale=.75,baseline=-3cm]
\draw (0,0) \br (-2,2) node[above] {$\lambda(a,b)^*$};
\draw (1,0) \br (2.75,2) node[above] {$\lambda(ab,c)^*$};
\begin{scope}[xshift=-2cm,yshift=-4.5cm]
  \foreach \x in {0,1}{
  \draw(\x,0)--(\x,2);
  }
  \draw[fill=white] (-.25,.75) rectangle node {$a,b,c$} (1.25,1.25);
\end{scope}
\draw (-2,-4.5) -- (-2,-10);
\draw (-1,-4.5)--(-1,-9);
\draw (0,-10) \br (-1,-9);
\draw (0,-10) to [out=-90, in=-90] (5,-10)--(5,-7);
\draw[looseness=2] (-2,-10) to [out=-90, in=-90] (-1,-10);
\draw[white, line width=10] (-1,-10) \br (0,-9) \br (4,-7);
\draw (-1,-10) \br (0,-9) \br (4,-7);
\draw (-2,-2.5) \br (-3.5,2) node[above] {$\lambda(a,b)$};
\draw[white, line width=10] (-1,-2.5) --(-1,-1.5) \br (1,2);
\draw (-1,-2.5)--(-1,-1.5) \br (1,2) node[above] {$\lambda(ab,c)$};
\draw[looseness=1., white, line width=10] (-3,-10) node[below] {$Z_c$} --(-3,-4) \br (-.5,2);
\draw[looseness=1] (-3,-10) node[below] {$Z_c$} --(-3,-4) \br (-.5,2) node[above] {$Z_c$}; 
\begin{scope}[scale=1,xshift=2cm,yshift=-4.5cm]
\draw (0,0)--(0,2);
    \draw (1,0)--(1,2);
    \draw[fill=white] (-.25,.75) rectangle (1.25,1.25);
    \draw (.5,1) node {\small $a,b,c$};
    \draw (1.25,1.25) node[right] {$^{-1}$};
    \draw[looseness=1.5] (0,2) to [out=90,in=90] (3,2) -- (3,-1);
    \draw[looseness=1.5] (1,2) to [out=90,in=90] (2,2) -- (2,-1);
    \draw[looseness=1.5] (0,0) to [out=-90,in=-90] (-1,0) -- (-1,3);
    \draw[looseness=1.5] (1,0) to [out=-90,in=-90] (-2,0) -- (-2,3);
    \draw (3,-2.5) \br (2,-1);
     \draw[white, line width=10]  (2,-2.5) \br (3,-1);
     \draw  (2,-2.5) \br (3,-1);
    \draw (-1,3) \br (-2,4.5);
    \draw[white, line width=10] (-2,3) \br (-1,4.5);
      \draw (-2,3) \br (-1,4.5);
\end{scope}
    \end{tikzpicture} = 
\begin{tikzpicture}[line width=1,baseline=2cm]
\draw[looseness=2] (2.5,4) node [above] {$\lambda(a,b)$} to [out=-90,in=-90] (3.75,4) node [above] {$\lambda(a,b)^*$};
\draw (5,0) node [below] {$Z_c$}--(5,4) node [above] {$Z_c$};
\draw[looseness=2] (6.25,4)node [above] {$\lambda(ab,c)$}  to [out=-90,in=-90] (7.5,4)node [above] {$\lambda(ab,c)^*$} ;
    \end{tikzpicture}.
\end{align}
After applying a Reidemeister II move to the left picture and using naturality the $\nu_{a,b,c}$ and its inverse isomorphism can be canceled. Further Reidemeister II moves then give the identification with the picture on the right. A similar calculation shows that the monoidal functor $(\Id_{\CC},J)$ is braided. Since the twists are preserved by the identity functor it is clear the braided tensor equivalence is ribbon.

Now we will call two zestings of $\CC$ are {\it equivalent} if they result in the same zested category. That is, $(\lambda,\nu,t,f) \simeq (\lambda',\nu',t',f')$ if $\CC^{\lambda,\nu,t,f} \simeq \CC^{\lambda',\nu',t',f'}$. 

\begin{rmk} Note that this notion of equivalence (when restricted to braided zesting data) is more general than the notion of similarity recalled earlier in Section \ref{sec:notationassumptions}.
\end{rmk}

It is natural to ask whether the product operation introduced above induces induces a group structure on equivalence classes of zesting data, although we leave that for future work. 

\section{Reshetikhin-Turaev invariants of framed links factorize under zesting}
\label{sec:zestinglinks}

Next we examine another nice property of ribbon zesting, namely that it transforms the Reshetikhin-Turaev invariants of framed links in a well-behaved manner. Let $\CC$ be a strict $A$-graded ribbon fusion category and $(\lambda, \nu, t,f)$ the data of a ribbon zesting. Then for a fixed framed link $L$, the Reshetikhin-Turaev invariant assigned to $L$ colored by simple objects in $\CC^{\lambda,\nu,t,f}$ factors as a product of the Reshetikhin-Turaev invariant of $L$ colored by the same objects in $\CC$ and a ``zested part".

The proof of this statement will be clear after describing how to compute the Reshetikhin-Turaev invariants of framed links in $\CC^{\lambda,\nu,t,f}$ (Section \ref{sec:zestinglinkdiagrams}) and working out some concrete examples (Examples \ref{examp:W} and \ref{examp:B}). Since the Whitehead link and Borromean rings play a role in the study of the Mignard-Schauenburg modular isotopes in \cite{BDGRTW, KMS} and will reappear in Section \ref{sec:zestingmodularisotopes}, it will be instructive to use them as the links in these examples. 

\subsection{Zesting link diagrams}
\label{sec:zestinglinkdiagrams}

Let $L$ be a link diagram presented as the Markov closure of a braid diagram. In order to be consistent with previous sections we take the convention that the Artin generators $\sigma_i$ of the $n$-strand braid group $B_n$ and their inverses are drawn 
\begin{align*}
\sigma_i \quad = \quad
\begin{tikzpicture}[line width=1,baseline=10]
\draw (-1,0) node[below] {$1$}--(-1,1) node[above] {$1$};
\draw (-.5,.5) node {$\cdots$};
\brrev
\draw (0,0) node[below] {$i+1$};
\draw (1,0) node[below] {$i$};
\draw (0,1) node[above] {$i$};
\draw (1,1) node[above] {$i+1$};
\draw (2,0) node[below] {$n$} --(2,1) node[above] {$n$};
\draw (1.5,.5) node {$\cdots$};
\end{tikzpicture} \qquad & \qquad \sigma_i^{-1} \quad = \quad \begin{tikzpicture}[line width=1,baseline=10]
\draw (-1,0) node[below] {$1$}--(-1,1) node[above] {$1$};
\draw (-.5,.5) node {$\cdots$};
\brpos
\draw (0,0) node[below] {$i+1$};
\draw (1,0) node[below] {$i$};
\draw (0,1) node[above] {$i$};
\draw (1,1) node[above] {$i+1$};
\draw (2,0) node[below] {$n$} --(2,1) node[above] {$n$};
\draw (1.5,.5) node {$\cdots$};
\end{tikzpicture}.
\end{align*}
 We read braid words $b=\sigma_{i_1}^{m_1}\sigma_{i_2}^{m_2}\cdots\sigma_{i_k}^{m_k}$ for $m_j \in \mathds{Z}$ from left to right and multiplication in $B_n$ is given by stacking going down the page.

Let $L$ be the Markov closure of the $n$-strand braid word $b=\sigma_{i_1}^{m_1}\sigma_{i_2}^{m_2}\cdots\sigma_{i_k}^{m_k}$ where the $i$th strand is colored by the simple object $X_i$ in a strict ribbon fusion category $\CC$.

Recall from Section \ref{sec:linkinvariants} that the Reshetikhin-Turaev invariant of $L$ is the trace of the isomorphism associated to its underlying braid representative $b$, which can be expressed in $\CC$ as

$$\rho = \overset{k}{\underset{j=1}{\circ}}   \left ( \id_{X_1 \otimes \cdots \otimes X_{i_j -1 }} \otimes c_{X_{i_j,X_{i_j+1}}}^{\text{sgn}(m_j)} \otimes \id_{X_{j_1 +2} \otimes \cdots \otimes X_{n}} \right ).$$

Now consider a ribbon zesting of $\CC$. The morphism $\rho^{\lambda,\nu,t,f}$ representing $b$ in the (no longer strict) category $\CC^{\lambda,\nu,t}$ can be depicted using string diagrams in $\CC$ and the diagrams encoding the zesting data from Section \ref{sec:zesting} as follows. 

Suppose $X \in \CC_a$, $Y \in \CC_b$, and $Z \in \CC_c$. The braid isomorphisms $c_{X,Y}$ get replaced with the isomorphisms $c^{\lambda,t}_{X,Y}$

\begin{align*}
c_{X,Y} = &\begin{tikzpicture}[line width=1,scale=.75, baseline=.75cm]
\begin{scope}[yscale=2]
\brrev
\draw (0,0) node[below] {$\parcen{Y}$};
\draw (1,1) node[above] {$\parcen{Y}$};
\draw (1,0) node[below] {$\parcen{X}$};
\draw (0,1) node[above] {$\parcen{X}$};
\end{scope}
\end{tikzpicture} \qquad \mapsto \qquad c_{X,Y}^{\lambda,t} = \begin{tikzpicture}[line width=1,scale=.75, baseline=.75cm]
\begin{scope}[yscale=2]
\brrev
\draw (0,0) node[below] {$\parcen{Y}$};
\draw (1,1) node[above] {$\parcen{Y}$};
\draw (1,0) node[below] {$\parcen{X}$};
\draw (0,1) node[above] {$\parcen{X}$};
\end{scope}
\begin{scope}[xshift=2.5cm]
\draw (0,0) node[below] {$\lambda(b,a)$} --(0,2) node[above] {$\lambda(a,b)$}; 
\draw[fill=white] (0,1) node {$\scriptstyle t(a,b)$} circle (.5) ;
\end{scope}
\end{tikzpicture} 
\end{align*}
and the previously trivial associators $\alpha_{X,Y,Z}$ are replaced with $\alpha^{\lambda,\nu}_{X,Y,Z}$. 
\begin{align*}
\alpha_{X,Y,Z} = \begin{tikzpicture}[line width=1,scale=.75,baseline=.5cm]
\draw (0,-.5)node[below] {$X$}--(0,2) node[above] {$X$};
\draw (1,-.5) node[below] {$Y$} \br (1,2) node[above] {$Y$};
\draw (2,-.5)node[below] {$Z$} --(2,2)node[above] {$Z$};
\end{tikzpicture}  \mapsto  & \quad \alpha^{\lambda,\nu}_{X,Y,Z} = 
\begin{tikzpicture}[scale=.75,baseline=.5cm, line width=1]
\zestedassoc{a}{b}{c}{0}
\draw (0,-.5) node [below] {$\parcen{X}$} -- (0,0);
\draw (1,-.5) node [below] {$\parcen{Y}$} -- (1,0);
\draw (2,-.5) node [below] {$\parcen{Z}$} -- (2,0);
\draw (3,-.5) node [below] {$\lambda(b,c)$} -- (3,0);
\draw (4,-.5)  -- (4,0);
\draw (4.5,-.5) node [below] {$\lambda(a,bc)$};
\draw (0,2) node [above] {$\parcen{X}$};
\draw (1,2) node [above] {$\parcen{Y}$};
\draw (2,2) node [above] {$\lambda(a,b)$};
\draw (3,2) node [above] {$\parcen{Z}$};
\draw (4,2) node [above] {$\lambda(ab,c)$};
\end{tikzpicture}.
\end{align*}
We call the resulting diagram of $\rho^{\lambda,\nu,t}$ the zested braid diagram of $b$. 

One further obtains a zested link diagram of $L$ by taking the zested trace of the zested braid diagram. From first principles, this can be achieved by replacing the cups and caps of $\CC$ with those in $\CC^{\lambda,\nu,t,f}$ 
\begin{align*}
\begin{tikzpicture}[scale=1, line width=1,baseline=-10]
\draw[looseness=1.5] (0,0) node[above] {$X_a^{*}$} to [out=-90, in=-90] (2,0) node[above] {$\phantom{\overline{X_a}}X_a\phantom{\overline{X_a}}$};
\end{tikzpicture}
& \mapsto \qquad  \begin{tikzpicture}[scale=1, line width=1,baseline=-10]
\draw (1,0) node[above] {$\lambda(a,-a)^{*}$};
\draw[looseness=1.5] (1,0)  to [out=-90, in=-90] (3,0);
\draw (3,0) node[above] {$\phantom{\lambda(a)^{*}}\lambda(-a,a)\phantom{\lambda(a)^{*}}$};
\draw[fill=white] (2.7,-.1) rectangle node {\scriptsize $\nu$} (3.1, -.5);
\draw (3.1, -.2) node[right]  {$^{-1}$};
\draw[looseness=1.5,white, line width=10] (0,0) to [out=-90, in=-90] (2,0);
\draw[looseness=1.5] (0,0) node[above] {$X_a^{*}$} to [out=-90, in=-90] (2,0) node[above] {$\phantom{\overline{X_a}}X_a\phantom{\overline{X_a}}$}; 
\end{tikzpicture}\\
\begin{tikzpicture}[scale=1,line width=1,baseline=0]
\draw[looseness=1.5] (0,0) node[below] {$\phantom{X_a^{*}}X_a\phantom{X_a^{*}}$} to [out=90, in=90] (2,0) node[below] {$X_a^{*}$};
\end{tikzpicture} & \mapsto \begin{tikzpicture}[scale=1.25,line width=1,baseline=0]
\draw[looseness=1.5] (0,0) node[below] {$\phantom{X_a^{*}}X_g\phantom{X_a^{*}}$} to [out=90, in=90] (1,0) node[below] {$X_a^{*}$};
\draw (1.875,0) node[below] {$\lambda(a,-a)^{*}$};
\draw[looseness=1.5] (2,0) to [out=90, in=90] (3,0)  ;
\draw (3.125,0) node[below] {$\phantom{\lambda(a)^{*}}\lambda(a,-a)\phantom{\lambda(a)^{*}}$};
\end{tikzpicture}
\end{align*}

and inserting the picture for the zested pivotal isomorphism from \cite[Equation 5.6]{DGPRZ}.

Equivalently however one can use the original trace in $\CC$ and correct by a factor of $f/t$, see Equation \ref{eq:zestedtracejtrivial}. In other words, one can use the cups, caps, and (trivial) pivotal structure in $\CC$ to obtain a diagram which is equal to the zested link diagram after multiplying by $\frac{f(\sum_i^n a_i)}{t(\sum_i^n a_i,\sum_i^n a_i)}$, where the simple $X_i$ labeling the $i$th strand lives in $\CC_{a_i}$. We will do the latter to save space in drawing the zested diagrams and to simplify graphical arguments.

\begin{examp}
\label{examp:W}
Consider the braid word $b = \sigma_1^{-2}\sigma_2\sigma_1^{-1}\sigma_2 \in B_3$ representing the Whitehead link $\widetilde{W}$ with blackboard framing. Labeling the components of $\widetilde{W}$ by simple objects $X_i \in \CC_i$, $Y_j \in \CC_j$, $i,j \in A$ induces a labeling of $b$ as shown on the left hand side below. Viewing $b$ as a morphism in $\rho \in \text{End}_{\CC}((X_i \otimes Y_j) \otimes X_i)$ we draw the zested morphism $ \rho^{\lambda,\nu,t} \in \text{End}_{\CC^{\lambda,\nu,t,f}}((X_i \overset{\lambda}{\otimes} Y_j) \overset{\lambda}{\otimes} X_i)$ by replacing the crossings and inserting associators as described in Section \ref{sec:zestinglinkdiagrams} above:

    \begin{align*}
   \begin{tikzpicture}[line width=1,scale=.65,baseline=5cm]
   \draw (1,0) node[below] {$Y_j$}  \br (0,1);
   \draw [white, line width=10](0,0) \br (1,1);
   \draw (0,0)  node[below] {$X_i$} \br (1,1);
   \begin{scope}[yshift=1cm]
   \draw (1,0) \br (0,1);
   \draw [white, line width=10](0,0) \br (1,1);
   \draw (0,0) \br (1,1);
   \end{scope}
   \draw (2,0) node[below] {$X_i$} --(2,4);
   \draw (0,2) --(0,7);
   \draw (1,2)--(1,4);
\begin{scope}[xshift=1cm,yshift=4cm]
\draw (0,0) \br (1,1);
\draw [white, line width=10](1,0) \br (0,1);
\draw (1,0) \br (0,1);
\end{scope}
   \begin{scope}[yshift=7cm]
   \draw (1,0) \br (0,1);
   \draw [white, line width=10](0,0) \br (1,1);
   \draw (0,0) \br (1,1);
   \end{scope}
\begin{scope}[xshift=1cm,yshift=10cm]
\draw (0,0) \br (1,1);
\draw [white, line width=10](1,0) \br (0,1);
\draw (1,0) \br (0,1);
\end{scope}
   \draw (0,8)--(0,13) node[above] {$\phantom{Y_j}X_i\phantom{Y_j}$};
   \draw (1,5)--(1,7);
   \draw (2,5)--(2,10);
   \draw (1,8)--(1,10);
   \draw (1,11)--(1,13) node[above] {$Y_j$};
   \draw (2,11)--(2,13) node[above] {$\phantom{Y_j}X_i\phantom{Y_j}$};
   \end{tikzpicture} & \mapsto  \begin{tikzpicture}[line width=1,scale=.65,baseline=5cm]
\draw (0,0) node[below] {$X_i$};
   \draw (1,0) node[below] {$Y_j$};
   \draw (2.125,0) node[below] {$\lambda(i,j)$};
   \draw (3.25,0) node[below] {$X_i$};
   \draw (4.5,0) node[below] {$\lambda(ij,i)$};
   \zestedbraidminv{0}{\scriptstyle i,j}
   \zestedbraidminv{1}{\scriptstyle j,i}
    \zestedassocinv{i}{j}{i}{2}
    \zestedbraid{4}{\scriptstyle i,j}
    \zestedassoc{i}{i}{j}{5}
    \zestedbraidminv{7}{\scriptstyle i,i}
    \zestedassocinv{i}{i}{j}{8}
    \zestedbraid{10}{\scriptstyle j,i}
    \zestedassoc{i}{j}{i}{11}
    \draw (0,13) node[above] {$\parcen{X_i}$};
    \draw (1,13) node[above] {$\parcen{Y_j}$};
     \draw (2.125,13) node[above] {$\lambda(i,j)$};
    \draw (3.25,13) node[above] {$\parcen{X_i}$};
    \draw (4.5,13) node[above] {$\lambda(ij,i)$};
    \end{tikzpicture}
    \end{align*}
   After applying the formula for the zested trace from Equation \ref{eq:zestedtracejtrivial} it follows from naturality and Reidemeister II that the trace of the zested part of the diagram is disjoint from the trace of $\rho$ and gives a morphism in $\Hom(\mathds{1},\mathds{1})$. The corresponding scalar can be evaluated by simply multiplying together all of the coupon labels in the picture above. In particular, all of the terms from the $\nu$ isomorphisms cancel.
   \begin{align}
\label{eq:W}
\widetilde{W}^{\lambda,\nu,t,f}_{X_i,Y_j}  = \Tr^{t,f}(\rho^{\lambda,\nu,t})  = \frac{f(2i+j)}{t(2i+j,2i+j)} \Tr(\rho^{\lambda,\nu,t}) = {f(2i+j)}{t(2i+j,2i+j)} t(i,i)^{-1} \widetilde{W}_{X_i,Y_j} \end{align}
     
We will be interested in the closely related $W$-matrix that comes from changing the framing of $\widetilde{W}$ so that $W_{X,Y} = \theta_Y^{-1}\widetilde{W}_{X,Y}$ \cite{BDGRTW}. The zested $W$-matrix is related to the original $W$-matrix by
\begin{align*}
W^{\lambda,\nu,t,f}_{X_i,Y_j} &= f(-j)(\theta_{Y_j}^{-1})\cdot \widetilde{W}^{\lambda,\nu,t,f}_{X_i,Y_j}
&= f(-j)(\theta_{Y_j}^{-1}) \frac{f(2i+j)}{t(2i+j,2i+j)} t(i,i)^{-1}\widetilde{W}_{X_i,Y_j}=\frac{f(-j)f(2i+j)}{t(2i+j,2i+j)t(i,i)} W_{X_i,Y_j}.\end{align*}
\end{examp}

\begin{examp}
\label{examp:B}
The Borromean rings can be represented by $b=(\sigma_1\sigma_2)^3 \in B_3$. Coloring the strands by $X_i \in \CC_i, Y_j \in \CC_j, Z_k \in \CC_k$ the braid diagram transforms as follows.

\begin{align*}
\begin{tikzpicture}[line width=1,scale=.65,baseline=.65*9cm]
\brpos
\draw (0,1)--(0,6);
\draw (1,1)--(1,3);
\draw (2,0)--(2,3);
\begin{scope}[xshift=1cm,yshift=3cm]
\brrev
\draw (0,1)--(0,3);
\draw (1,1)--(1,3);
\end{scope}
\begin{scope}[yshift=6cm]
\brpos
\draw (0,1)--(0,6);
\draw (1,1)--(1,3);
\draw (2,0)--(2,3);
\begin{scope}[xshift=1cm,yshift=3cm]
\brrev
\draw (0,1)--(0,3);
\draw (1,1)--(1,3);
\end{scope}
\end{scope}
\begin{scope}[yshift=12cm]
\brpos
\draw (0,1)--(0,6) node[above] {$\phantom{_j}X_i\phantom{_j}$};
\draw (1,1)--(1,3);
\draw (2,0)--(2,3);
\begin{scope}[xshift=1cm,yshift=3cm]
\brrev
\draw (0,1)--(0,3) node[above] {$Y_j$};
\draw (1,1)--(1,3) node[above] {$\phantom{_j}Z_k\phantom{_j}$};
\end{scope}
\end{scope}
\draw (0,0) node [below] {$\phantom{_j}X_i\phantom{_j}$};
\draw (1,0) node [below] {$Y_j$};
\draw (2,0) node [below] {$\phantom{_j}Z_k\phantom{_j}$};
\end{tikzpicture}
& \quad \mapsto \quad  \begin{tikzpicture}[line width=1,scale=.65,baseline=.65*9cm]
    \zestedbraidminv{0}{\scriptstyle i,j}
    \zestedassocinv{j}{i}{k}{1}
    \zestedbraid{3}{\scriptstyle k,i}
    \zestedassoc{j}{k}{i}{4}
    \zestedbraidminv{6}{\scriptstyle j,k}
    \zestedassocinv{k}{j}{i}{7}
    \zestedbraid{9}{\scriptstyle i,j}
    \zestedassoc{k}{i}{j}{10}
    \zestedbraidminv{12}{\scriptstyle k,i}
    \zestedassocinv{i}{k}{j}{13}
    \zestedbraid{15}{\scriptstyle j,k}
    \zestedassoc{i}{j}{k}{16}
    \draw (0,18) node[above] {$\phantom{_j}X_i\phantom{_j}$};
    \draw (1,18) node[above] {$Y_j$};
    \draw (2.125,18) node[above] {$\lambda(i,j)$};
     \draw (3.25,18) node[above] {$\phantom{_j}Z_k\phantom{_j} $};
    \draw (4.5,18) node[above] {$\lambda(ij,k)$};
    \draw (0,0) node [below] {$\phantom{_j}X_i\phantom{_j}$};
\draw (1,0) node [below] {$Y_j$};
\draw (2.125,0) node [below] {$\lambda(i,j)$};
\draw (3.25,0) node [below] {$\phantom{_j}Z_k\phantom{_j}$};
\draw (4.5,0)  node [below] {$\lambda(ij,k)$};
\end{tikzpicture}
\end{align*}

This time we find that the contributions from the $t$ isomorphisms cancel and we are left with
    \begin{align}
    \label{eq:B}
        B^{\lambda,\nu,t,f}_{X_i,Y_j,Z_k} & = \Tr^{t,f}(\rho^{\lambda,\nu,t})  = \frac{f(i+j+k)}{t(i+j+k,i+j+k)} \Tr(\rho^{\lambda,\nu,t})= \frac{f(i+j+k)}{t(i+j+k,i+j+k)}\frac{\nu_{i,j,k}\nu_{k,i,j}\nu_{j,k,i}}{\nu_{i,k,j}\nu_{k,j,i}\nu_{j,i,k}}B_{X_i,Y_j,Z_k}.
    \end{align}
\end{examp}

The kind of formulas derived in Examples \ref{examp:W} and \ref{examp:B} are much simpler for cyclic zestings, i.e. when $A \cong \mathds{Z}/N\mathds{Z}$. In this case, each of $\nu$, $t$, and $f$ can be expressed in terms of a single parameter $s$ as recorded in Section \ref{sec:cyclictannakianzesting}. 
\begin{prop}[Formulas for some topological invariants under cyclic zesting]
\label{prop:cycliczestedinvariants}
Let $\CC$ be a ribbon fusion category and $(\lambda_a,\nu_b,t_s,f_s)$ a cyclic ribbon zesting of $\CC$.
\begin{enumerate}
    \item $W^{\lambda_a,\nu_b,t_s,f_s}_{X_i,Y_j}=s^{i^2+j^2} W_{X_i,Y_j}$
    \item $B^{\lambda_a,\nu_b,t_s,f_s}_{X_i,Y_j,Z_k}= B_{X_i,Y_j,Z_k}$
    \end{enumerate}
for all $X_i \in \CC_i, Y_j \in \CC_j, Z_k \in \CC_k$.
\end{prop}
\begin{proof}
The formulas follow from Equations \ref{eq:W} and \ref{eq:B} using that $\nu(i,j,k)$ is symmetric in $i$ and $j$ for cyclic zestings and substituting $t(i,j)=s^{-ij}$ and $f(i)=s^{-i^2}$, see Equations \ref{eq:3cochain}-\ref{eq:cycliczestingf} in Section \ref{sec:cyclictannakianzesting}.
\end{proof}

This invariance of the Borromean tensor $B_{X,Y,Z}$ under cyclic zesting is an important fact that we will need later in the proof of Theorem \ref{thm:zestingmodularisotopes}. 
\begin{rmk}
\label{rmk:zesting3manifoldinvariants1}
Recall that when $\CC$ is modular the Reshetikhin-Turaev construction also gives well-defined invariants of closed, oriented $3$-manifolds $M$. Given a presentation of a fixed $M$ by surgery on a Kirby-colored framed link $L$, one can compute the zested 3-manifold invariant essentially by summing up the contributions from the zested link invariants colored by simple objects in $\CC$ according to the formula in \cite[Equation 2.2.a]{TuraevBook}.
\end{rmk}
The following example shows that in principle it is straightforward to compute specific zested $3$-manifold invariants.
\begin{examp}
\label{ex:lens}
Let $\CC$ be a modular fusion category. We will further assume that $\CC$ is unitary, so that in particular all of the quantum dimensions are positive. The Lens spaces $L_{p,q}$ where $(p,q)=1$ can be obtained by surgery on the $n$-component link 
\begin{align*}
\begin{tikzpicture}[baseline=0,scale=.5]
\begin{scope}[thick,scale=.55]
\draw (0,0) arc (0:-30:1);
\draw (0,0) arc (0:300:1);
\draw (1.375,0) arc (0:120:1);
\draw (1.375,0) arc (0:-30:1);
\draw (.875,-.875) arc (-60:-210:1);
\draw (1.75,1) arc (90:120:1);
\draw (.875,.5) arc (150:270:1);
\draw (2.5,0) node {\large$\cdots$};
\draw (4.25,0) arc (0:90:1);
\draw (4.25,0) arc (0:-30:1);
\draw (3.75,-.875) arc (-60:-90:1);
\draw (5.625,0) arc (0:120:1);
\draw (5.625,0) arc (0:-210:1);
\end{scope}
\end{tikzpicture}
\end{align*}
 whose $i$th component is framed by $a_i$, where $n$ is the length of the continued fraction representation $$
-p / q=a_n-\frac{1}{a_{n-1}-\frac{1}{\ddots \frac{1}{a_2-\frac{1}{a_1}}}}
.$$
The Reshetikhin-Turaev invariant of $L(p,q)$ is computed by the formula
\begin{align}
\label{eq:lens}
Z_{\CC}\left(L(p, q)\right)=\frac{\Theta^{-\sigma(L)}}{D^{n+1}} \sum_{X_1, \ldots, X_n}\left(\prod_{j=1}^n d_{X_j}\left(\theta_{X_j}\right)^{a_j}\right) \frac{S_{X_1 X_2}}{d_{X_2}} \frac{S_{X_2 X_3}}{d_{X_3}} \ldots \frac{S_{X_{n-2} X_{n-1}}}{d_{X_{n-1}}} S_{X_{n-1} X_n} 
\end{align}
where $\Theta=\frac{1}{D} \sum\limits_{X\in \Irr(\CC)} d_X^2 \theta_X$ is the central charge and $\sigma(L)$ is the signature of the matrix $$
\left(\begin{array}{ccccccc}
a_1 & -1 & 0 & \ldots & 0 & 0 & 0 \\
-1 & a_2 & -1 & \ldots & 0 & 0 & 0 \\
0 & -1 & a_3 & -1 & 0 & \ldots & 0 \\
\ldots & \ldots & \ldots & \ldots & \ldots & \ldots & \ldots \\
\ldots & \ldots & \ldots & \ldots & \ldots & a_{n-1} & -1 \\
0 & 0 & \ldots & 0 & 0 & -1 & a_n
\end{array}\right),
$$
see \cite[Section 2.2.4]{BDGRTW}. 

Note that $Z_{\CC}(L_{p,q})$ depends only on the modular data of $\CC$. Thus given zesting data $(\lambda,\nu,t,f)$ for $\CC$ a formula for $Z_{\CC^{\lambda,\nu,t,f}}(L_{p,q})$ follows easily from plugging in the formulas for the zested modular data from Equations \ref{eq:zestedmodulardata} into Equation \ref{eq:lens} above. Since $\CC$ is unitary, $\dim(\lambda(a,b))=1$ for all $a,b \in A$ and the global quantum dimension $D$ is preserved. This simplifies the formulas for the zested modular data considerably.

\begin{align*} Z_{\CC^{\lambda,\nu,t,f}}(L_{p,q}) =\frac{(\Theta^f)^{-\sigma(L)}}{D^{n+1}} \sum_{X_1,\ldots, X_n}\left(\prod_{j=1}^n d_{X_j}(f(b_j)\theta_{X_j})^{a_j}\right )\left (\prod_{k=1}^{n-1}\frac{t(b_k,-b_{k+1})^2f(b_k-b_{k+1})S_{X_k,\lambda(b_k,-b_{k+1})}S_{X_k,X_{k+1}}}{t(b_k-b_{k+1},b_k-b_{k+1})d_{X_{k+1}}}\right )d_{X_n}
\end{align*}
where $X_j \in \CC_{b_j}$, $\Theta^f = \frac{1}{D} \sum_{X_b \in \Irr(\CC)} d_{X_b}^2f(b)\theta_{X_b}$, and $t(b_k,-b_{k+1})^2= t(b_k,-b_{k+1})t(-b_{k+1},b_k)$.
\end{examp}
\subsection{Zesting Reshetikhin-Turaev invariants of framed links}
The method illustrated in the examples of the previous section can be generalized to analyze how an arbitrary Reshetikhin-Turaev link invariant transforms under zesting. 

\begin{thm}
\label{thm:zestinglinks} Let $\CC$ be an $A$-graded ribbon fusion category, let $\mathcal{C}^{\lambda,\nu,t,f}$ be a ribbon zesting of $\CC$, and let $L$ be a framed $n$-component link whose $j$-th component is colored by $X_{i_j} \in \Irr(\CC_{i_j})=\Irr(\mathcal{C}^{\lambda,\nu,t,f}_{i_j})$. Then
\begin{align}
    L^{\lambda,\nu,t,f}_{X_{i_1}, X_{i_2}, \ldots, X_{i_n}}&=J(\lambda, \nu, t,f, i_1, i_2, \ldots, i_n)\cdot L_{X_{i_1}, X_{i_2}, \ldots, X_{i_n}}
\end{align}
where $J$ is a scalar-valued function of the ribbon zesting data $(\lambda, \nu, t,f)$ and the gradings $i_j$.
\end{thm}

\begin{proof}
Let $L$ be an $N$-component link and let $b \in B_n$ be a braidword representing $L$. Fixing a coloring of the $N$ components of $L$ by simple objects in an $A$-graded ribbon category $\CC$ induces a labeling of the strands of $b$ interpreted as a morphism in $\CC$, say $X_{i_1}, X_{i_2}, \ldots, X_{i_n}$, where $X_{i_j} \in \CC_{i_j}$. By abuse of notation we will continue to conflate the Reshetikhin-Turaev invariant associated to $L$ with $L$.

Now let $(\lambda, \nu, t,f)$ be the data of a ribbon zesting of $\CC$. Then the braid isomorphism corresponding to $b$ transforms under the zesting according to the diagram below, where we have suppressed the dependence on $\nu$ and $t$.

\begin{center}
\begin{tikzpicture}[line width=1,yscale=.5, xscale=1.125]
\draw (0,0) node[below] {
\small $\parcen{X_{i_1}}$} --(0,12) node[above] {\small $\parcen{X_{i_1}}$};
\draw (1,0) node[below] {\small $\parcen{X_{i_2}}$} --(1,12) node[above] {\small $\parcen{X_{i_2}}$};
\draw (2,0) node[below] {\small$ (i_1,i_2)$}  \br  (5,4) -- (5,8)  \br (2,12) node[above] {\small$ (i_1,i_2)$};
\draw (4,0)  \br (6,4) --(6,8) \br (4,12);
\draw (4,0) node[below] {\small $ (i_1+i_2,i_3)$};
\draw (4,12) node[above] {\small $ (i_1+i_2,i_3)$};
\draw[white, line width=5] (3,0) \br (2,4) -- (2,8) \br (3,12) ;
\draw (3,0) node[below] {\small$\parcen{X_{i_3}}$} \br (2,4) -- (2,8) \br (3,12) node[above] {\small$\parcen{X_{i_3}}$};
\draw[white, line width=10] (5,0) \br (3,4)--(3,8) \br (5,12);
\draw (5,0) node[below] {\small $\parcen{X_{i_4}}$} \br (3,4)--(3,8) \br (5,12) node[above] {\small $\parcen{X_{i_4}}$};
\draw (6.5,0) node[below] {\scriptsize $\cdots \, (\sum_{j=1}^{n-2} i_j,i_{n-1})$};
\draw (7,0) \br (8,4) --(8,8) \br (7,12);
\draw (6.5,12) node[above] {\scriptsize $ \cdots \,(\sum_{j=1}^{n-2} i_j,i_{n-1})$};
\draw[white, line width=10] (8,0) \br (4,4)-- (4,8) \br (8,12) ;
\draw (8,0) node[below] {\small $\parcen{X_{i_n}}$} \br (4,4) --(4,8) \br (8,12) node[above] {\small $\parcen{X_{i_n}}$};
\draw[fill=white] (-.5,5) rectangle node {$b$} (4.5,7); 
\draw (3.5,4.5) node {\Large $\cdots$};
\draw (3.5,7.5) node {\Large $\cdots$};
\draw (9.5,0) node[below] {\tiny$ (\sum_{j=1}^{n-1} i_j,i_n)$};
\draw (9,0) --(9,12);
\draw (9.5,12) node[above] {\tiny $ (\sum^{n-1}_{j=1} i_j,i_n)$};
\draw (5.875,1) node[below] {\Large $\cdots$};
\draw (5.875,11) node[above] {\Large $\cdots$};
\end{tikzpicture}
\end{center}

The key observation is that the part of the zested isomorphism on the strands labeled by $\lambda$ lie under the original braid isomorphism. As explained in Section \ref{sec:zestinglinkdiagrams}, taking the zested categorical trace of this picture gives the zested link invariant, and this zested trace is equal to the original trace by Equation \ref{eq:zestedtracejtrivial}, so that up to the scalar factor $f(\sum_{j=1}^N i_j)/t(\sum_{j=1}^N i_j,\sum_{j=1}^N i_j)$ the zested invariant is given by the Markov closure of the picture in dashed red lines below. Moreover, by sphericality we can draw the Markov closure so that the tracial closure of the part of the diagram involving strands labeled by $\lambda$ lies completely under the Markov closure of $b$.

\begin{tikzpicture}[line width=1,yscale=.25, xscale=.5,baseline=30]
\foreach \x in {0,1,2,3,5}
{\draw (\x,0) --(\x,12);}
\draw[fill=white] (-.5,5) rectangle node {$b$} (4.5,7); 
\draw (4,2.5) node {\Large $\cdots$};
\draw (4,9.5) node {\Large $\cdots$};
\draw[looseness=2] (0,12) to [out=90,in=90] (-1,12)--(-1,0) to [out=-90,in=-90] (0,0);
\draw[looseness=2] (1,12) to [out=90,in=90] (-1.5,12)--(-1.5,0) to [out=-90,in=-90] (1,0);
\draw[looseness=2] (2,12) to [out=90,in=90] (-2,12)--(-2,0) to [out=-90,in=-90] (2,0);
\draw[looseness=2] (3,12) to [out=90,in=90] (-2.5,12)--(-2.5,0) to [out=-90,in=-90] (3,0);
\draw[looseness=2] (5,12) to [out=90,in=90] (-3.5,12)--(-3.5,0) to [out=-90,in=-90] (5,0);
\end{tikzpicture}\qquad \Huge $\mapsto$ \normalsize
\begin{tikzpicture}[line width=1,yscale=.25, xscale=.5,baseline=30]
\draw (0,0) --(0,12);
\draw (1,0)  --(1,12);
\draw (2,0) \br  (5,4) -- (5,8)  \br (2,12);
\draw (4,0)  \br (6,4) --(6,8) \br (4,12);
\draw[white, line width=5] (3,0) \br (2,4) -- (2,8) \br (3,12) ;
\draw (3,0) \br (2,4) -- (2,8) \br (3,12);
\draw[white, line width=10] (5,0) \br (3,4)--(3,8) \br (5,12);
\draw (5,0)  \br (3,4)--(3,8) \br (5,12);
\draw (7,0) \br (8,4) --(8,8) \br (7,12) ;
\draw[white, line width=10] (8,0) \br (4,4)-- (4,8) \br (8,12) ;
\draw (8,0)\br (4,4) --(4,8) \br (8,12);
\draw[fill=white] (-.5,5) rectangle node {$b$} (4.5,7); 
\draw (3.5,4.5) node {\Large $\cdots$};
\draw (3.5,7.5) node {\Large $\cdots$};
\draw (9,0)  --(9,12);
\draw (5.875,0) node[below] {\Large $\cdots$};
\draw (5.875,12) node[above] {\Large $\cdots$};
\draw[red, dashed] (-.75,0) rectangle (9.25,12);
\draw[looseness=2] (0,12) to [out=90,in=90] (-1,12)--(-1,0) to [out=-90,in=-90] (0,0);
\draw[looseness=2] (1,12) to [out=90,in=90] (-1.5,12)--(-1.5,0) to [out=-90,in=-90] (1,0);
\draw[looseness=2] (2,12) to [out=90,in=90] (11,12)--(11,0) to [out=-90,in=-90] (2,0);
\draw[white, line width=10, looseness=2] (3,12) to [out=90,in=90] (-2,12)--(-2,0) to [out=-90,in=-90] (3,0);
\draw[looseness=2] (3,12) to [out=90,in=90] (-2,12)--(-2,0) to [out=-90,in=-90] (3,0);
\draw[looseness=2] (9,12) to [out=90,in=90] (9.5,12)--(9.5,0) to [out=-90,in=-90] (9,0);
\draw[looseness=2] (7,12) to [out=90,in=90] (10,12)--(10,0) to [out=-90,in=-90] (7,0);
\draw[looseness=2] (4,12) to [out=90,in=90] (10.5,12)--(10.5,0) to [out=-90,in=-90] (4,0);
\draw[white, line width=10, looseness=2] (5,12) to [out=90,in=90] (-2.5,12)--(-2.5,0) to [out=-90,in=-90] (5,0);
\draw[looseness=2] (5,12) to [out=90,in=90] (-2.5,12)--(-2.5,0) to [out=-90,in=-90] (5,0);
\draw[white, line width=10, looseness=2] (8,12) to [out=90,in=90] (-3,12)--(-3,0) to [out=-90,in=-90] (8,0);
\draw[looseness=2] (8,12) to [out=90,in=90] (-3,12)--(-3,0) to [out=-90,in=-90] (8,0);
\end{tikzpicture}

Now observe that the zested link diagram can be isotoped so into two disjoint diagrams, the Markov closure of $b$ and a ``zested part" involving only strands labeled by $\lambda$. In equation form this says that 

$$L^{\lambda,\nu, t,f}_{X_{i_1}, X_{i_2}, \ldots, X_{i_N}}=J(\lambda, \nu, t,f, i_1, i_2, \ldots, i_N)\cdot L_{X_{i_1}, X_{i_2}, \ldots, X_{i_N}}$$
where $J(\lambda, t,f, i_1, i_2, \ldots, i_N)$ is the scalar obtained from evaluating the ``under" diagram and multiplying by the aforementioned factor of $f(\sum_{j=1}^N i_j)/t(\sum_{j=1}^N i_j,\sum_{j=1}^N i_j)$.
\end{proof}

\begin{rmk}
\label{rmk:zesting3manifoldinvariants2}
When $\CC$ is modular one cannot conclude that the 3-manifold invariants factorize in the same manner as in Theorem \ref{thm:zestinglinks} due to their definition involving labeling the components of the surgery diagram by the Kirby color. 
\end{rmk}

\section{Zesting produces modular isotopes}
\label{sec:zestingmodularisotopes}

In this section, we show that the Mignard-Schauenburg modular isotopes can be realized by the ribbon zesting construction. These modular isotopes belong to a family of modular fusion categories arising from finite groups and 3-cocycles which can be understood in two equivalent ways. First, as Drinfeld centers $\mathcal{Z}(\Vect_G^{\omega})$ of $G$-graded vector spaces with associativity constraint $\omega \in Z^3(G,\mathds{C}^{\times})$, and second as the categories $\Rep(D^{\omega}G)$ of finite dimensional modules over the twisted quantum double Hopf algebras $D^{\omega}G$.\footnote{Actually one must invert the 3-cocycle $\omega$ to make the identification, so that $\ZZ(\Vect_G^{\omega^{-1}}) \simeq \Rep(D^{\omega}G)$ \cite{GM}.} Both realizations of these categories can be convenient for understanding their properties and we will switch freely between them in this section. In either case we will denote them by $\mathcal{Z}(\Vect_G^{\omega})$ and abuse terminology by referring to the categories as twisted Drinfeld doubles of finite groups.

Section \ref{sec:drinfelddoubles} recalls formulas for the modular data of $\mathcal{Z}(\Vect_G^{\omega})$, while Section \ref{sec:msmodularisotopes} establishes some facts specific to the Mignard-Schauenburg modular isotopes which occur at particular values of $G$ and $\omega$. In Section \ref{sec:zestingproducesmodularisotopes} we give explicit ribbon zesting data that constructs the Mignard-Schauenburg modular isotopes using the understanding of how Reshetikhin-Turaev invariants of framed links transform under zesting developed in Section \ref{sec:zestinglinks}. We conclude in Section \ref{sec:conclusions} with some further remarks and open questions.

\subsection{Modular data for twisted Drinfeld doubles of finite groups}
\label{sec:drinfelddoubles}

Isomorphism classes of simple objects in $\mathcal{Z}(\Vect^{\omega}_G)$ are parametrized by pairs $([t], \rho_t)$ where $[t]$ is a conjugacy class in $G$ with a fixed representative $t$ and $\rho_t$ is an irreducible $\theta_t$-projective representation of the centralizer $C_G(t) = \{g \in G | gtg^{-1}=t\}$ with factor set
\begin{equation}
    \theta_t(x,y) = \frac{\omega(t,x,y)\omega(x,y,(xy)^{-1}txy)}{\omega(x,x^{-1}tx,y)}.
\end{equation}
Write $\chi_t$ for the character of the representation $\rho_t$.
The quantum dimensions of simple objects are given by
\begin{align}
\label{qdims}
     d_{([t], \rho_t)}= |[t] |\dim(\rho_t)
\end{align}
and the modular data by
\begin{align}
\label{eq:twisteddoublesmat}
 S_{([a], \rho_a),([b],\psi_b)} &= 
    \frac{1}{|G|} \sum_{\substack{g \in [a]\\ h \in [b] \cap C_G(g)}} \chi^*_{\rho_g}(h)\chi^*_{\psi_h}(g)\\
    \label{eq:twisteddoublesmat2}
      \begin{split}  &= \frac{1}{|G|} \sum_{\substack{g \in [a]\\ h \in [b] \cap C_G(g)}} \left (\frac{\theta_a(a_g,h)\theta_a(a_gh,a_g^{-1})\theta_b(b_h,g)\theta_b(b_hg,b_h^{-1})}{\theta_g(a_g^{-1},a_g)\theta_h(b_h^{-1},b_h)} \right)^*  \\ & \qquad \times \chi^*_{\rho_a}(a_gha_g^{-1})\chi^*_{\psi_b}(b_hgb_h^{-1}) \end{split}\\
     \label{eq:twisteddoubletmat}
    T_{([a],\rho_a),([b],\psi_b)} &= \delta_{[a],[b]}\delta_{\rho_a,\psi_b}\frac{\chi_{\rho_a(a)}}{\chi_{\rho_a(e)}},
\end{align}
  where $g=a_gaa_g^{-1}$ and $h=b_hbb_h^{-1}$.
The first formula for the $S$-matrix in in Equation \ref{eq:twisteddoublesmat}  is more concise \cite{CGR}, but Equation \ref{eq:twisteddoublesmat2} has the benefit of being readily computable from the data of the simple objects (a priori one has $\psi_a$ for a fixed conjugacy class representative $a$ but not necessarily $\psi_g$ for all $g \in [a]$ \cite{GM}).  

\subsection{The Mignard-Schauenburg modular isotopes}
\label{sec:msmodularisotopes}

The Mignard-Schauenburg categories are of the form $\ZZ(\Vect_G^{\omega})$ where $G=\mathds{Z}_q \rtimes \mathds{Z}_p$ for odd primes $p>3$ and $q$ such that $p | q-1$. There are $|H^3(G,\mathds{C}^{\times})|=p$ such categories up to equivalence, each with rank $\frac{q^2-1+p^3}{p}$ and Frobenius-Perron dimension $p^2q^2$ (global quantum dimension $pq$). From now on let $\omega$ denote a normalized $3$-cocycle whose cohomology class generates $ H^3(G,\mathds{C}^{\times})\cong \mathds{Z}/p\mathds{Z}$ and fix representatives $\omega^u$ of the $p$ different cohomology classes of $G$ where $u \in \{0,1,\ldots, p-1\}$. 

While each of the $p$ cohomology classes corresponds to a distinct modular fusion category $\mathcal{Z}(\Vect_G^{\omega^u})$, there are only 3 distinct sets of modular data up to a relabeling \cite{MS}. One set of modular data corresponds to $u=0$. The other two sets of modular data are each shared by $(p-1)/2$ different modular fusion categories, which we call modular isotopes. We will need further technical details about the Mignard-Schauenburg modular data (see Section \ref{sec:msmodulardata}), but first we recall the definition of modular isotopes and what is known about how to distinguish them with invariants beyond the modular data.

Recall that two sets of modular data $(S^{(1)},T^{(1)})$ and $(S^{(2)}, T^{(2)})$ are equal up to a relabeling if there is a permutation $\pi$ of isomorphism classes of simple objects whose associated permutation matrix $P$ maps the two sets of modular data onto each other, i.e.~
\begin{align*}
S^{(1)} = PS^{(2)}P^{-1}\\
T^{(1)} = PT^{(2)}P^{-1}.
\end{align*}
In other words, the $S$- and $T$-matrices are related by a permutation change-of-basis matrix. When two sets of modular data are equal up to a relabeling one cannot rule out that $\pi$ lifts to a ribbon tensor autoequivalence of $\CC^{(1)}$ and $\CC^{(2)}$ from the modular data alone. This is the case with each of the two families of $\frac{p-1}{2}$ Mignard-Schauenburg modular isotopes: their $S$- and $T$-matrices aren't equal with respect to the natural ordered basis of $\Irr(\mathcal{Z}(\Vect_G^{\omega^{u}}))$ introduced in Section \ref{sec:msmodulardata}, but there exist permutations of $\Irr(\mathcal{Z}(\Vect_G^{\omega^{u}}))$ for certain values of $u$ which can be applied to identify their modular data.

Therefore, exhibiting modular isotopes requires an additional invariant beyond the modular data that obstructs a lifting of a relabeling $\pi$ to a ribbon equivalence. In other words, one must demonstrate an invariant that fails to respect all possible relabelings that induce identifications of the modular data. 

Let $B$ be the Borromean tensor of invariants associated to a modular fusion category from Example \ref{examp:B} and $W$ the matrix of invariants associated to the Whitehead link from Example \ref{examp:W}.

\begin{prop}[Theorem 6.2 in \cite{KMS}]
\label{prop:WB5}
The Mignard-Schauenburg modular isotopes are distinguished by
$(B,T)$ for any $p,q$.
\end{prop} \qed
\begin{examp}[\cite{BDGRTW}]
The $p=5, q=11$ Mignard-Schauenburg isotopes are distinguished by $(W,T)$.
\end{examp}

\begin{rmk}
A subtle point is that the Borromean tensor is actually identical for each choice of $u$ with respect to the ordered basis of $\text{Irr}(\ZZ(\Vect^{\omega^u}_G))$ \cite{KMS}, unlike $S$,$T$, or $W$. Nevertheless it is not preserved by any of the relabelings that identify the twists, so each $(B,T)$ still corresponds uniquely to one of the modular isotopes.  
\end{rmk}

Of course, one is not limited to using knot and link invariants for the purpose of distinguishing modular isotopes, see \cite{WW19} where they use traces of mapping class group representations of higher genus surfaces. 

\subsubsection{General properties of Mignard-Schauenburg categories}

The next two propositions record some facts we will need about these categories. The first will be used to establish that $\ZZ(\Vect^{\omega^u}_G)$ has the prerequisites to apply the zesting construction, namely a grading with nontrivial invertible objects in its trivial component.

\begin{lem}
\label{prop:MSpointed}
Let $\mathcal{C}=\ZZ(\Vect^{\omega^u}_G)$ where $G=\mathds{Z}_q \rtimes \mathds{Z}_p$, with $p>3$ and $q$ odd primes such that $p|q-1$. Then
\begin{enumerate} 
    \item $\rank(\mathcal{C}_{pt})=p$,
    \item $\mathcal{C}_{pt} \subset \mathcal{C}_{ad}$, and
    \item $\mathcal{C}_{pt}$ is Tannakian.
\end{enumerate}
\end{lem}

\begin{proof}
(1) By Equation \ref{qdims}, an object $([t], \rho_t)\in \mathcal{C}=\ZZ(\Vect^{\omega^u}_G)$ is invertible if and only if $|[t]| = 1$ and $\dim(\rho_t) = 1$. But this is equivalent to having $t\in Z(G) = \langle e \rangle$ and $\dim(\rho_t) = 1$. In this case, we have $\theta_t = \theta_e = 1$, since $\omega^u$ can be chosen to be a normalized $3$-cocycle. Then all projective representations of $C_G(e) = G$ are linear representations. Since the number of one-dimensional representations of $G$ is $|G/[G,G]| = p$, the claim holds. 

(2) It follows from (1) that $\FPdim (\mathcal{C}_{ad}) = pq^2$. By the Burnside Theorem for fusion categories (\cite[Theorem 1.2]{ENO2}), $\mathcal{C}_{ad}$ is solvable and $(\mathcal{C}_{ad})_{pt}$ is nontrivial, by \cite[Proposition 4.5 (iii)]{ENO2}. Since  $(\mathcal{C}_{ad})_{pt}$ is a subcategory of $\CC_{pt}$ with $\FPdim \CC_{pt} = p$, then $(\mathcal{C}_{ad})_{pt} = \CC_{pt}$ as desired.  

(3) The pointed subcategory $\mathcal{C}_{pt}$ is symmetric by (2) and  \cite[Corollary 6.9]{GN}. Moreover, it follows from (1) that $\FPdim(\mathcal{C}_{pt}) = p$, which is an odd prime, and therefore $\mathcal{C}_{pt}$ is Tannakian by \cite[Corollary 2.50]{DGNO}. 
\end{proof}

The next proposition says that the Mignard-Schauenburg categories are (up to equivalence) the unique modular fusion categorifications of their fusion rules.

\begin{prop}
\label{prop:msmodulardata}
If a modular fusion category $\mathcal B$ has the same Grothendieck ring as a Mignard-Schauenburg category $\CC = \mathcal Z(\Vect_G^{\omega})$, with $G=\mathds{Z}_q \rtimes \mathds{Z}_p$ and $p$, $q$  odd primes such that $p | q-1$, then $\mathcal B$ is in fact a Mignard-Schauenburg category, that is $\mathcal B  \simeq \mathcal Z(\Vect_G^{\tilde{\omega}})$ for some $\tilde{\omega}$.
\end{prop}

\begin{proof} 
The first step of the proof is to show that $\mathcal B$ is the Drinfeld center of some fusion category. To show this, it is enough to find a connected \'etale algebra in $\mathcal B$ such that the square of its Frobenius-Perron dimension equals $\FPdim (\mathcal B)$, see \cite[Corollary 4.1 (i)]{DMNO}. We will find such an \'etale algebra by looking for a copy of $\Rep(G)$ in $\mathcal{B}$ and constructing its function algebra \cite[Example 2.8]{DMNO}.

Notice that since $\mathcal B$ has the same Grothendieck ring as $\CC \simeq \Rep(D^{\omega}(G))$ and $\Rep (G)$ is a fusion subcategory of $\CC$, then there exists a fusion subcategory $\mathcal D$ of $\mathcal B$ of dimension $pq = |G|$. By Proposition \ref{prop:MSpointed}(i), $\mathcal D$ cannot be pointed, and, since $p$ and $q$ are odd, it follows from \cite[Theorem 6.3(i)]{EGO} that $\mathcal D \cong \Rep(G)$ as fusion categories. Moreover, since $\mathcal B$ is braided so is $\mathcal D$. 
The possible braidings on $\Rep(G)$ were described on \cite{Dav_qt}, \cite[Example 5.5]{Nik_braidings}; they are in bijection with triples $(L, M, B)$, where $L$ and $M$ are
a pair of mutually commuting normal subgroups of $G$ and $B:L\times M \to \textbf{k}^{\times}$ is a non-degenerate
$G$-invariant bilinear form. For $G=\mathds{Z}_q \rtimes \mathds{Z}_p$, the possible choices for $L$ and $M$ are for both to be the trivial subgroup $\{e\}$ or both $\mathds{Z}_q$ but this second choice is not possible under the condition of admitting a $G$-invariant bilinear form. So $M = L = \{e\}$ which implies that the braiding on $\mathcal D \simeq \Rep(G)$ is symmetric.

In this case, we have $\mathcal B \simeq \mathcal Z(\mathcal D)$ for a fusion category $\mathcal D$ of Frobenius-Perron dimension $pq$. 
It follows from the classification of fusion categories of Frobenius-Perron dimension $pq$ given in \cite[Theorem 6.3]{EGO} that $\mathcal D \simeq \Vect_{\tilde{G}}^{\tilde{\omega}}$ for some $\tilde{G}$ with $|\tilde{G}|=pq$ and 3-cocycle $\tilde{\omega}$.

Suppose the group $\tilde{G}$ were isomorphic to $\mathbb Z_{pq}$. The Schur multiplier of  $\mathbb Z_{pq}$ is the trivial Schur multiplier because all its Sylow subgroups are cyclic. So $\mathcal Z(\Vect_{\mathbb Z_{pq}})$ is pointed and therefore cannot be equivalent to $\mathcal B$. Thus $\mathcal B \simeq Z(\Vect_G^{\tilde{\omega}})$ where $G \cong \mathds{Z}_p \rtimes \mathds{Z}_p$ as desired.
\end{proof}

\subsubsection{Mignard-Schauenburg modular data}
\label{sec:msmodulardata}
We now give a more detailed description of the modular data for $\ZZ(\Vect^{\omega^u}_{\mathds{Z}_q \rtimes \mathds{Z}_p})$. To express the semidirect product structure explicitly we take $G=\mathds{Z}_q \rtimes_n \mathds{Z}_p$ where $n^p \equiv 1 \mod q$ and $n\ne 1$. Elements are written $g=a^lb^k$
where $\langle a \rangle =\mathds{Z}_q$ and $\langle b \rangle = \mathds{Z}_p$, and multiplication is given by

\begin{align*}a^lb^k\cdot a^{l'}b^{k'}= a^{l+n^kl'}b^{k+k'}.\end{align*}

Let $g=a^{l_g}b^{i_g}$ and $[x]_p=x \mod p$. The 3-cocycles $\omega^u$ take the following form, see Equation (4.6) of \cite{BDGRTW}.

\begin{align}
\label{eq:ms3cocycle}
    \omega^u(g,h,k)&= e^{\frac{2\pi i}{p^2}u[i_k]_p([i_g]_p+[i_h]_p-[i_g+i_h]_p)}\\
    &= \begin{cases}1 & \text{if } [i_g]_p+[i_h]_p < p \\ e^{\frac{2\pi i}{p}u[i_k]_p} & \text{if } [i_g]_p+[i_h]_p \ge p \end{cases}.
\end{align}
Recall that simple objects in $\ZZ(\Vect^{\omega^u}_G)$ then correspond to pairs $([t],\rho_t)$ where $\rho_t$ is a $\theta^u_t$-projective representation of $C_G(t)$.  Since the $3$-cocycle does not depend on the power of $a$ in its arguments, it follows that $\theta^u_t$ is also independent of $a$, and we may write

\begin{equation}
\label{eq:2cochain}
   \theta^u_{b^k}(b^{i_1},b^{i_2}):= \theta^u_{a^lb^k}(a^{l_1}b^{i_1},a^{l_2}b^{i_2}) = \omega^u(b^{i_1},b^{i_2},b^k)  =  \begin{cases}1 & \text{if } i_1 + i_2 < p \\ e^{\frac{2\pi i}{p}uk} & \text{if } i_1 + i_2 \ge p \end{cases}
\end{equation}
where $i_1, i_2 \in \{0,1,\ldots, p-1\}$

Recall from Lemma \ref{prop:MSpointed}(1) that there are $p$ invertible simple objects. Since $\mathcal{Z}(\Vect_G^{\omega^u})$ is modular, it follows that its universal grading group is isomorphic to $\mathds{Z}/p\mathds{Z}$ \cite[Lemma 8.22.9]{EGNO}. The invertibles -- corresponding to the $p$ linear one-dimensional representations of $G$ -- all lie in the trivial component of the universal grading by Lemma \ref{prop:MSpointed}(2). The universal grading of an arbitrary simple object $([t],\rho_t)$ can read off from its conjugacy class according to the following table.

\begin{table}[H]
        \label{eq:conjugacyclasses}
    \begin{tabular}{c|c|c}
    Type of Conjugacy Class & Centralizer & Universal $\mathds{Z}/p\mathds{Z}$-Graded Component \\
\hline
    $[e]$ & $G$ & 0 \\

    $[a^l]$, \quad $l \in \{1,2, \ldots, q-1\}$ & $\mathds{Z}_q$ & 0 \\

  $[b^k]$,\quad $k \in \{1, 2, \ldots, p-1\}$ & $\mathds{Z}_p$ & $k$ \\
    \end{tabular}
    \end{table}

Since the objects with $[t]=[b^k]$ for some $k\ne 0$ correspond to the nontrivial components of the universal grading they play an important role in our analysis. In particular, we will use that the $\theta^u_t$-projective irreps of their centralizers $\mathds{Z}/p\mathds{Z}$ are given explicitly by
\begin{align}
\label{eq:projirreps}
\pi^{u,s}_k(b^l)=\underbrace{e^{\frac{2\pi i}{p}sl}}_{{\text{linear part}}}\cdot \underbrace{e^{\frac{2\pi i}{p^2}ukl}}_{\text{projective part}},
\end{align}
see \cite[Equation 4.7]{BDGRTW}. In what follows we will implicitly use the lexicographic order in $k$ and $s$ respectively to fix an order of the objects $([b^k], \pi^{u,s}_k)$. \\

We will need the following technical lemma establishing the structure of the block of the $S$-matrix coming from the double braiding of the simple objects of the form $([b^k],\pi^{u,s}_k)$. 

\begin{lem}
\label{lem:MSSmatformula}
The nontrivially-graded blocks of the $S$-matrix for the Mignard-Schauenburg modular isotopes satisfy
\begin{align} S_{([b^{k_1}], \pi^{u,s}_{k_1}),([b^{k_2}],\pi^{u,s}_{k_2})} = \frac{1}{|G|} \sum_{\substack{g \in [b^{k_1}]\\ h \in [b^{k_2}] \cap C(b^{k_1})}}  \pi^{u, s_1}_{k_1}(a_gha_g^{-1})^*\pi^{u, s_2}_{k_2}(b_hgb_h^{-1})^*.\end{align}
Moreover, the normalization factor involving the $\theta$'s in Equation \ref{eq:twisteddoublesmat2} can be safely ignored for all simple objects.
\end{lem}

\begin{proof}
Note that since the representations are 1-dimensional we can replace the characters with the projective irreps from Equation \ref{eq:projirreps} themselves. \\
Consider two simple objects in the $k_1$- and $k_2$-graded component of $\mathcal{Z}(\Vect^{\omega^u}_{\mathds{Z}_q \rtimes_n \mathds{Z}_p})$, where $k_1,k_2 \in \{1, 2, \ldots, p-1\}$. Such objects are the form $([b^{k_1}],\pi_{k_1}^{u,s_1})$ and $([b^{k_2}],\pi_{k_2}^{s_2})$ respectively. Then by the discussion preceding Equation \ref{eq:2cochain} we have
\begin{align*}
    \frac{\theta^u_a(a_g,h)\theta^u_a(a_gh,a_g^{-1})\theta^u_b(b_h,g)\theta^u_b(b_hg,b_h^{-1})}{\theta^u_g(a_g^{-1},a_g)\theta^u_h(b_h^{-1},b_h)} \\
   = \frac{\theta^u_{b^{k_1}}(b^{k_1'},b^{k_2})\theta_{b^{k_1}}(b^{[k_1'+k_2]_p},b^p-k_1')\theta_{b^{k_2}}(b_{k_2'},b^{k_1})\theta_{b^{k_2}}(b^{[k_2'+k_1]_p},b^{p-k_2'})}{\theta_{b^{k_1}}(b^{p-k_1'},b^{k_1'})\theta_{b^{k_2}}(b^{p-k_2'},b^{k_2'})} \\
\end{align*}
where we have identified $a=b^{k_1},b=b^{k_2}$. Remembering only the power of $b$ in the elements $g, h, a_g$, and $ b_g$, where $g=a_gaa_g^{-1}$ and $h=b_hbb_h^{-1}$, we identify $g=b^{k_1}, h=b^{k_2}$ and put $a_g=b^{k_1'}, b_h=b^{k_2'}$ for some $k_1', k_2'$. Since the formula for $\theta^u_{b^k}$ in Equation \ref{eq:2cochain} requires its arguments to take values in $\mathds{Z}/p\mathds{Z}$ it is necessary to take the residue mod $p$ when the powers of $b$ might exceed $p$. \\
Next observe that $\theta^u_{b^{k_1}}(b^{k_1'},b^{k_2})\theta_{b^{k_1}}(b^{[k_1'+k_2]_p},b^p-k_1') = e^{\frac{2\pi i}{p}uk_1}$. If $k_1'+k_2 <p$ then $[k_1'+k_2]_p+p-k_1'=k_1'+k_2+p-k_1'=k_2+p \ge p$, and the second factor contributes a factor of $e^{\frac{2\pi i}{p}uk_1}$. If $k_1'+k_2 \ge p$ then the first factor contributes a factor of $e^{\frac{2\pi i}{p}uk_1}$ and the second factor is trivial since $[k_1'+k_2]_p+p-k_1'=(k_1'+ k_2 - p)+p -k_1' = k_2 <p$. An identical argument for $\theta_{b^{k_2}}(b_{k_2'},b^{k_1})\theta_{b^{k_2}}(b^{[k_2'+k_1]_p},b^{p-k_2'})$ shows that the numerator will always contribute a factor of $e^{\frac{2\pi i}{p}uk_1}e^{\frac{2\pi i}{p}uk_2}$, which then gives
\begin{align*} \frac{\theta^u_{b^{k_1}}(b^{k_1'},b^{k_2})\theta_{b^{k_1}}(b^{[k_1'+k_2]_p},b^p-k_1')\theta_{b^{k_2}}(b_{k_2'},b^{k_1})\theta_{b^{k_2}}(b^{[k_2'+k_1]_p},b^{p-k_2'})}{\theta_{b^{k_1}}(b^{p-k_1'},b^{k_1'})\theta_{b^{k_2}}(b^{p-k_2'},b^{k_2'})} 
=\frac{e^{\frac{2\pi i}{p}uk_1}e^{\frac{2\pi i}{p}uk_2}}{e^{\frac{2\pi i}{p}uk_1}e^{\frac{2\pi i}{p}uk_2}} 
  =1
\end{align*}
as desired. One can also verify that ratio of $\theta^u$'s vanishes for the simple objects from conjugacy classes $[e]$ and $[a^l]$.
\end{proof}

\subsection{Constructing the Mignard-Schauenburg categories via zesting}
\label{sec:zestingproducesmodularisotopes}
Next we show that when the ribbon zesting construction is applied to $\mathcal{Z}(\Vect_G)$, it produces $p$ categories whose modular data is equal to the modular data of the Mignard-Schauenburg categories. 
\begin{prop}
\label{prop:mszesting}
Let $\CC^{(u)}= \ZZ(\Vect_G^{\omega^u})$ denote the Mignard-Schauenburg categories and $S^{(u)}$, $T^{(u)}$ their corresponding $S$- and $T$-matrices. Then for all $u \in \mathds{Z}/p\mathds{Z}$ there is a choice of cyclic zesting data $(\lambda_a,\nu_b,t_s,f_s)$ for which

\begin{align*} (S^{(0)})^{\lambda_a,\nu_b,t_s,f_s} &= S^{(u)}\\
(T^{(0)})^{\lambda_a,\nu_b,t_s,f_s} &= T^{(u)}
\end{align*}
In other words, the $p$ modular fusion categories obtained from zesting $\CC^{(0)}$ have the same modular data as the $p$ categories $\CC^{(u)}$.  
\end{prop}

\begin{proof}
Since $\mathcal{Z}(\Vect^{\omega^u}_G)$ is modular, the universal grading group is $\mathds{Z}/p\mathds{Z}$. By Lemma \ref{prop:MSpointed}, the pointed subcategory of $\mathcal{Z}(\Vect^{\omega^u}_G)$ is cyclic Tannakian, i.e. of the form $\Rep(\mathds{Z}/p\mathds{Z})$, and is contained in $\CC_{ad}$, the trivially component of the universal grading. Therefore we can apply the theory of cyclic Tannakian zestings with respect to the universal $\mathds{Z}/p\mathds{Z}$-grading from Section \ref{sec:cyclictannakianzesting}. 

The choices of associative zestings of $\mathcal{Z}(\Vect_G)$ correspond to pairs $(a,b) \in \mathds{Z}/p\mathds{Z} \times\mathds{Z}/p\mathds{Z}$ satisfying $a \equiv 2b \mod p$. Since $p$ is prime for each $a$ there is a unique $b$ for which the obstruction to associative zesting vanishes, and each of these admits exactly one ribbon zesting. We may then consider the family of ribbon zestings of the form $(\lambda_{2u}, \nu_u,t_s, f_s)$ controlled by the parameter $u \in \{1,\ldots,p-1\}$. Here $s=e^{-2\pi i u/p^2}$
where
$$s^p=\zeta^{-2u}, \quad \text{where $\zeta$ is a $2p$th-root of unity}.$$
Recall from Equations \ref{eq:cyclictannakianSmatrix} and \ref{eq:cyclictannakianTmatrix} of Section \ref{sec:cyclictannakianzesting} how the modular data transforms under cyclic Tannakian zesting:

\begin{align*}
S_{X_i,Y_j}^{\lambda_{2u}, \nu_u,t_s, f_s}=s^{2ij}S_{X_i,Y_j}, &&  0\leq i ,j < N,  X_i\in \CC_i, Y_j\in \CC_j.
\end{align*}
\begin{align*}
T_{X_i,X_i}^{\lambda_{2u}, \nu_u,t_s, f_s}= s^{-i^2}T_{X_i,X_i}, &&  0\leq i < N, X_i\in \CC_i,
\end{align*}

We must show that the modular data arising from ``twisting" of the untwisted Drinfeld double $\mathcal{Z}(\Vect_G)$ changes in precisely this same way.

According to Equation \ref{eq:ms3cocycle}, changing the power of the $3$-cocycle $\omega^u$ only affects the twists of anyons of the form $([b^k], \chi)$ according to the formula                                             
\begin{align}
    \theta^{(u)}_{([b^k],\chi)} &= e^{2\pi i k^2 u/p^2}  \theta_{([b^k],\chi)}. 
    \end{align}

The anyons $([b^k], \pi^{u,s}_{k})$ for $k=1,\, 2,\, \ldots, p-1$ are precisely the anyons in the $k$-graded components of the $\mathds{Z}/p\mathds{Z}$, which according to Equation \ref{eq:cyclictannakianTmatrix} change by a factor of $s^{-k^2}$ where we have identified the notation $i$ for indexing the graded components in \cite{DGPRZ} with the $k$ in \cite{MS} and \cite{BDGRTW}. Thus the way the twists transform under cyclic Tannakian zesting agrees with those of $\ZZ (\Vect^{\omega^u}_G)$.

To see that the $S$-matrices of the $\ZZ (\Vect^{\omega^u}_G)$ can be realized via ribbon zesting the $S$-matrix of $\ZZ\left (\Vect_G\right)$, we compare the formulas for the $S$-matrix of a quantum double in the twisted and untwisted case universally graded block by block.

It is straightforward to check that the entries of the twisted and untwisted $S$-matrices agree whenever one of the simples is in the trivially graded component. It then suffices to consider the nontrivially-graded blocks of the $S$-matrix, i.e. those entries corresponding to anyons of the form $([b^k], \pi^{u,s}_{k})$ for $k=1,\, 2,\, \ldots, p-1$. 

By Lemma \ref{lem:MSSmatformula}, these $S$-matrix entries obey the formula 

$$ S^{(u)}_{([b^{k_1}], \pi_{k_1}^{u,s_1}),([b^{k_2}], \pi_{k_2}^{u,s_2})} = \frac{1}{|G|} \sum_{\substack{g \in [b^{k_1}]\\ h \in [b^{k_2}] \cap C_G(g)}}  \pi_{k_1}^{u,s_1}(a_gha_g^{-1})^*\pi_{k_2}^{u,s_2}(b_hgb_h^{-1})^*$$

We can further simplify this formula using the explicit description of conjugacy classes in $\mathds{Z}_q \rtimes_n \mathds{Z}_p$ as in Equation \ref{eq:conjugacyclasses}, namely that $[b^{k_{1,2}}]=\{b^{k_{1,2}}, ab^{k_{1,2}}, a^2b^{k_{1,2}}, \ldots, a^{q-1}b^{k_{1,2}}\}$. Since $C_G(b^{k_{1,2}})=\langle b \rangle$, we must have $a_gha_g^{-1}=b^{k_2}$ and $b_hgb_h^{-1}=b^{k_1}$. In order words, the only conjugates of $h$ and $g$ that live in the centralizer of $b^{k_1}$ and $b^{k_2}$ respectively are $b^{k_2}$ and $b^{k_1}$ respectively. This gives

\begin{align*} S^{(u)}_{([b^{k_1}], \pi_{k_1}^{u,s_1}),([b^{k_2}], \pi_{k_2}^{u,s_2})} &= \frac{1}{|G|} \sum_{\substack{g \in [b^{k_1}]\\ h \in [b^{k_2}] \cap C_G(g)}}  \pi_{k_1}^{u,s_1}(a_gha_g^{-1})^*\pi_{k_2}^{u,s_2}(b_hgb_h^{-1})^*\\
& =\frac{1}{|G|} \sum_{\substack{g \in [b^{k_1}]\\ h \in [b^{k_2}] \cap C_G(g)}}  \pi_{k_1}^{u,s_1}(b^{k_2})^*\pi_{k_2}^{u, s_2}(b^{k_1})^*\\
& = e^{-\frac{2\pi i}{p^2}uk_1k_2}e^{-\frac{2\pi i}{p^2}uk_1k_2} \frac{1}{|G|} \sum_{\substack{g \in [b^{k_1}]\\ h \in [b^{k_2}] \cap C_G(g)}}  \pi_{k_1}^{s_1}(a_gha_g^{-1})^*\pi_{k_2}^{s_2}(b_hgb_h^{-1})^*\\
&=e^{-\frac{4\pi i}{p^2}uk_1k_2} S_{([b^{k_1}], \pi_{k_1}^{s_1}),([b^{k_2}], \pi_{k_2}^{s_2})}
\end{align*}
where we have used the explicit formulas of Equation \ref{eq:projirreps} for the projective characters as product of a linear character and projective part. These $S$-matrix entries for $\mathcal{Z}(\Vect^{\omega^u}_G)$ then differ by this overall factor of $e^{-\frac{4\pi i}{p^2}uk_1k_2}$ from the untwisted $S$-matrix for $\mathcal{Z}(\Vect_G)$. In summary, away from the trivially-graded anyons the $S$-matrix of the twisted Drinfeld double of $\mathds{Z}_q \rtimes_n \mathds{Z}_p$ differs from the untwisted double according to 

$$S^{(u)}_{([b^{k_1}], \pi_{k_1}^{u,s_1}),([b^{k_2}], \pi_{k_2}^{u,s_2})}=e^{-\frac{4\pi i}{p^2}uk_1k_2} S_{([b^{k_1}], \pi_{k_1}^{s_1}),([b^{k_2}], \pi_{k_2}^{s_2})}.$$

This factor of $e^{-\frac{4\pi i}{p^2}uk_1k_2}$ agrees exactly with the factor of $s^{2k_1k_2}=e^{-\frac{4\pi i}{p^2}uk_1k_2}$ given by zesting $\mathcal{Z}(\Vect_{\mathds{Z}_q \rtimes_n \mathds{Z}_p})$.
\end{proof}

\begin{thm}
\label{thm:zestingmodularisotopes}
The Mignard-Schauenburg modular isotopes are related by cyclic Tannakian zesting with respect to their universal grading group $\mathds{Z}/p\mathds{Z}$.
\end{thm}

\begin{proof}
By Proposition \ref{prop:mszesting}, the modular data of the Mignard-Schauenburg categories $\mathcal{Z}(\Vect^{\omega^u}_G)$ can be generated from zesting $\mathcal{Z}(\Vect_G)$. By Proposition \ref{prop:msmodulardata} we may conclude that these zestings of $\mathcal{Z}(\Vect_G)$ are in fact Mignard-Schauenburg categories. However, we know that modular isotopes occur among these categories, so one must rule out the possiblity that the $p$ different zestings only realizes a proper subset of the the $p$ different categories. Note though that by  Proposition \ref{prop:cycliczestedinvariants}(2) the $p$ cyclic zestings all share the same $B$-tensor as the the Mignard-Schauenburg categories. That is, in the notation of Proposition \ref{prop:mszesting} $(B^{(0)})^{\lambda_u,\nu_u, t_s,f_s} = B^{(u)}$ for all $u \in \mathds{Z}/p\mathds{Z}$. By Proposition \ref{prop:WB5}, \cite[Theorem 6.2]{KMS}, the $T$-matrix together with the $B$-tensor are enough to distinguish the different modular isotopes. Thus we may conclude that all $p$ Mignard-Schauenburg categories are realized by zesting $\mathcal{Z}(\Vect_G)$. 

It follows that any pair $\mathcal{Z}(\Vect_G^{\omega^u})$ and $\mathcal{Z}(\Vect_G^{\omega^{u'}})$ are also related by zesting: by the results in Section \ref{sec:zestingproperties} the product zesting $$\left ( \lambda_{u'}\cdot\lambda_u^*, \nu_{u'}\cdot\tilde{\nu}_u,t_{u'}\cdot (t_u^{-1})^*, f_{u'}\cdot(f_u^{-1})^* \right )$$ produces a category ribbon equivalent to $\mathcal{Z}(\Vect_G^{\omega^{u'}})$ from zesting $\mathcal{Z}(\Vect_G^{\omega^u})$. In particular, any pair of modular isotopes are related by zesting.
\end{proof}

This shows that zesting can produce different modular fusion categories with the same modular data. 

\begin{examp} The smallest Mignard-Schauenburg categories occur at rank 49 when $p=5,\, q=11$. We give the explicit choices of cyclic Tannakian zesting data for $\ZZ(\Vect_{\mathds{Z}_{11} \rtimes \mathds{Z}_5})$ that produce the 4 twisted doubles  $\ZZ(\Vect^{\omega^u}_{\mathds{Z}_{11} \rtimes \mathds{Z}_5})$, $u \in \{1,2,3,4\}$. 
\begin{center}
\begin{tabular}{c|c}
   \text{Associative zesting parameters} $(a,b)$ & Braided/ribbon zesting parameter $s$  \\
   \hline
   $(2,1)$  & $e^{-2\pi i/25}$ \\ 
   $(4,2)$  & $e^{-4\pi i/25}$ \\ 
   $(1,3)$  & $e^{-6\pi i/25}$ \\ 
   $(3,4)$  & $e^{-8\pi i/25}$ \\ 
\end{tabular}
\end{center}
One can check by computer that the modular data are related by $$S_{X_j,Y_k}^{(u)} = e^{-2\pi i ujk/25} \, S_{X_j, Y_k}$$
$$T_{X_j,X_j}^{(u)} = e^{2\pi i uj^2/25}\, T_{X_j, X_j}$$
while the $W$ and $B$ invariants satisfy
$$W_{X_j,Y_k}^{(u)} = e^{-2\pi i u(j^2+k^2)/25} W_{X_j, Y_k}$$
$$B^{(u)}_{X_j,X_k,X_l}=B_{X_j,X_k,X_l}.$$
\end{examp}

\subsection{Discussion}
\label{sec:conclusions}

Our proof of Theorem \ref{thm:zestingmodularisotopes} relies on specific properties of the Mignard-Schauenburg modular isotopes and explicit formulas for projective characters of certain groups. In general it will not be the case that an associative zesting of fusion categories taking $\Vect_G$ to $\Vect_G^{\omega}$ results in a braided zesting between the centers $\ZZ(\Vect_G)$ and $\ZZ(\Vect^{\omega}_G)$, as the following example shows.

\begin{examp}[See Example 3.5 of \cite{Ng2003}]
\label{examp:cesar}
Let $n >2$ and $q$ a primitive $n$th root of unity. Consider the pointed fusion category $\Vect_A^{\omega}$ with
\begin{align*}
A= \mathds{Z}/n\mathds{Z} \times \mathds{Z}/n\mathds{Z}\times \mathds{Z}/n\mathds{Z}, \quad\omega(\vec{a},\vec{b},\vec{c})=q^{\det(M)}
\end{align*}
where $\vec{a},\vec{b},\vec{c} \in A$ and $M = \begin{bmatrix} \vec{a} & \vec{b} & \vec{c} \end{bmatrix}$ is the $3\times 3$ matrix with columns given by $\vec{a},\vec{b},\vec{c}$. Then $\ZZ(\Vect_A)$ is pointed but $\ZZ(\Vect_A^\omega)$ is not. Since a zesting of a pointed category is again pointed it is not possible for $\ZZ(\Vect_A^\omega)$ to be a zesting of $\ZZ(\Vect_A)$.
\end{examp}

\begin{quest}
When does an associative zesting of a fusion category induce a braided zesting of its center?
\end{quest}

Theorem \ref{thm:zestingmodularisotopes} gives a new construction of the Mignard-Schauenburg modular isotopes and suggests that zesting may be a useful tool in the classification of modular fusion categories. 

It is natural to ask,
\begin{quest}
Do all modular isotopes arise from zesting?
\end{quest}

\begin{quest}
What is the smallest rank where modular isotopes occur?
\end{quest}
 
While Theorem \ref{thm:zestinglinks} allows one to predict the values of topological invariants under zesting, a more thorough understanding of when  different zestings produce distinct categories would be needed in order to distinguish the Mignard-Schauenburg modular isotopes directly from the zesting data without relying on additional invariants beyond the modular data.

If one could show formally that the different choices of zesting data used in Section \ref{sec:zestingmodularisotopes} necessarily lead to different categories then Theorem \ref{thm:zestingmodularisotopes} would follow immediately from Proposition \ref{prop:mszesting}.

\begin{quest}
\label{quest:whenzestinginequiv}
When do two choices of zesting data $(\lambda,\nu,t,f)$ and $(\lambda',\nu',t',f')$ produce equivalent categories?
\end{quest}

It is clear for example that when $\lambda$ and $\lambda'$ differ by a coboundary, the zested categories have isomorphic fusion rings, but the converse is not true. Take for example the eight unitary modular fusion categories with Ising fusion rules, some of which are related with nontrivial choices of $\lambda$. One must be similarly careful when addressing this question with the remainder of the zesting data.

\bibliographystyle{plain}
\bibliography{refs}
\end{document}